\documentclass[11pt,a4paper]{article}
\usepackage[english]{babel}
\usepackage{amsfonts}
\usepackage{amssymb}
\usepackage{epsfig}
\usepackage{graphicx}
\usepackage{amsmath}
\usepackage{amsthm}
\usepackage{enumitem}
\usepackage{tikz}
\usetikzlibrary{intersections}
\usetikzlibrary{arrows,arrows.meta,bending}

\usepackage{xcolor}
\definecolor{candyapplered}{rgb}{1.0, 0.03, 0.0}
\definecolor{coolblack}{rgb}{0.0, 0.18, 0.39}
\definecolor{darkpowderblue}{rgb}{0.0, 0.2, 0.6}
\definecolor{lightbrown}{rgb}{0.71, 0.4, 0.11}
\definecolor{ruddybrown}{rgb}{0.73, 0.4, 0.16}
\usepackage{hyperref}
\hypersetup{hyperfigures = true, colorlinks = true, linkcolor=darkpowderblue, citecolor=ruddybrown}

\usepackage[utf8]{inputenc}

\def\RR{\mathbb R}

\def\bet{\beta}

\def\del{\Delta}
\def\eps{\varepsilon}

\def\lam{\lambda}

\def\tet{\theta}

\newcommand{\pr}{\medskip\noindent\textit{\textbf{Proof.}} }

\newcommand{\prr}[1]{\medskip\noindent\textit{\textbf{Proof of Proposition #1.}} }

\newcommand{\prrr}[1]{\medskip\noindent\textit{\textbf{Proof of Theorem #1.}} }
\newcommand{\cqd}{{\unskip\nobreak\hfil\penalty50
		\hskip2em\hbox{}\nobreak\hfil\mbox{\rule{1ex}{1ex} \qquad}
		\parfillskip=0pt \finalhyphendemerits=0\par\medskip}}
	
\newtheorem{theo}{Theorem}[section]
\newtheorem{prop}{Proposition}[section]
\newtheorem{lem}{Lemma}[section]
\newtheorem{rem}{Remark}[section]

\date{}
\begin{document}
\title{Perturbation analysis in a free boundary problem arising in tumor growth model}
\author{Ahlem ABDELOUAHAB$^{a}$ and Sabri BENSID$^{b}$ }
\maketitle
\begin{center}
$^{a}$Dynamical Systems and Applications Laboratory\\
Department of Mathematics, Faculty of Sciences\\
University of Tlemcen, B.P. 119\\
Tlemcen 13000, Algeria\\
Mail: $ ahlemabdelouahab98@gmail.com$

$^{b}$Dynamical Systems and Applications Laboratory\\
Department of Mathematics, Faculty of Sciences\\
University of Tlemcen, B.P. 119\\
Tlemcen 13000, Algeria\\
Mail: $ edp\_sabri@yahoo.fr$
\end{center}
\begin{abstract}
We study the existence and multiplicity of solutions of the following free boundary problem
$$
(P)\left\{
\begin{array}{rcll}
\del u &=& \lam ( \eps +(1-\eps ) H(u-\mu))~ \hspace{3mm}&\text{in}~\Omega (t)\\
u&=& \overline{u}_{\infty}~\hspace{3mm} &\text{on } ~ \partial \Omega(t)
\end{array}
\right.
$$
where $\Omega(t) \subset \RR^3$ a regular domain at $t>0$,  $\eps,~\overline{u}_{\infty},~\lambda,~\mu$ are a positive parameters and $H$ is the Heaviside step function.
\\The problem (P) has two free boundaries: the outer boundary of  $\Omega(t)$ and the inner boundary  whose evolution is implicit generated by the discontinuous nonlinearity $H$. The  problem (P) arise in tumor growth models as well as in other contexts such as climatology. First, we show the existence and multiplicity of radial solutions of problem (P) where $\Omega(t)$ is a spherical domain. Moreover, the bifurcation diagrams are giving. Secondly, using the perturbation technic combining to local methods, we prove the existence of solutions and characterize the free boundaries of problem (P) near the corresponding radial solutions.
\end{abstract}

\section{Introduction}
Let us denote by $\Omega(t) \subset \RR^3 $  the spherical tumor region at $t >0$ and $R(t)$ its radius. The classical diffusion process (See \cite{Byrne Medicine}) shows that the oxygen (nutrient) concentration $u(x,t)$ satisfies the following equation, for $t>0$ and $x \in \Omega(t)$
$$\frac{\partial u}{\partial t}=D \del u -\lam f(u)$$
where $D >0$ is the diffusion coefficient and $\lam$ is a positive constant where $\lam f(u)$ describe the consumed oxygen rate.  There is a large theoretical literature devoted to mathematical models of the growth tumor with different forms of the nonlinearity $f$. In \cite{greenspan}, H.P. Greenspan was probably the first who proposed an avascular tumor growth with new depend variables for the local cell velocity and pressure within the tumor. He notices that the oxygen diffusion time scale is much shorter than a typical tumor doubling time to conclude that $u$ is quasi stable and satisfies the following problem, for $t >0$ and $x \in \Omega(t)$,
\begin{equation}\label{p1}
\left\{
\begin{array}{rcll}
D \del u &=& \lam f(u)~ \hspace{3mm}&\text{in}~\Omega (t)\\
u&=& \overline{u}_{\infty}~\hspace{3mm} &\text{on } ~ \partial \Omega(t),
\end{array}
\right.
\end{equation}
where $\overline{u}_{\infty}$ is the nutrient concentration at the boundary.
\\Many works have since extended Greenspan models to other growth scenario, see \cite{adam}, \cite{chaplain}, \cite{maggelakis} and \cite{mo} and the references therein. For example, McElwain and Morris \cite{mo} proposed the following nonlinearity
\begin{equation*}
f(u)=\hspace*{0.5cm}
\left\{
\begin{array}{rcll}
1&,& ~ \hspace{3mm}&\text{for}~u \geq \sigma_2\\
\frac{u}{\sigma_2}&,& ~\hspace{3mm} &\text{for}~ \sigma_1 < u \leq \sigma_2\\
0&,& ~\hspace{3mm} &\text{for}~ u \leq \sigma_1
\end{array}
\right.
\end{equation*}
where $\sigma_1$ and $\sigma_2$ are a positive critical values for which tumor goes from one phase to another. Byrne and Chaplain \cite{B} supposed that the consumed oxygen rate $f$ is given by $$f(u)=H(u-c)$$
where $H$ is the Heaviside step function $H(s) = 0 ~\text{for}~ s < 0;~ H(s) = 1 ~\text{for} ~s \geq 0$
\\and $c$ denotes the threshold oxygen concentration at which cells.
\\Noting that during the years 1995-1999, H. Byrne and M. Chaplain give a remarkable series of works when they analysed theoretically, numerically and by asymptotic analysis the above models. We refer the reader to \cite{adam}, \cite{B}, \cite{chaplain} and \cite{maggelakis} for more details.\\
\hspace*{3mm}Later, many authors developed and studied the evolution of tumor region in the form of free boundary problems and different models are explored. In particular the papers of A. Friedman et al. See \cite{cui friedman}, \cite{FR}, \cite{FR2} and the references given
there.\\\\
In this paper,  we are concerned to study the radial symmetric growth of problem (\ref{p1}). More precisely, taking the radial part of Laplace operator in three dimension, $D \equiv 1$ to simplified, we consider the following problem
\begin{equation}\label{p2}
\left\{
\begin{array}{rcll}
\vspace*{3mm} \frac{1}{r^2} \frac{\partial}{\partial r} \big ( r^2 \frac{\partial u}{\partial r}\big)&=& \lam f(u) ~ \hspace{3mm}&~0< r <R(t),\hspace*{2mm} t >0\\
\vspace*{3mm} u(R(t),t)&=&u_{\infty} ~\hspace*{3mm}&~ \text{for} \hspace*{2mm} t >0\\
\frac{\partial u}{\partial r}(0,t)&=&0~\hspace*{3mm}&~ \text{for} \hspace*{2mm} t >0,
\end{array}
\right.
\end{equation}
where \begin{eqnarray*}
f(u)&=& \eps +(1-\eps ) H(u-\mu)\\
&=&\left\{
\begin{array}{rcll}
\eps&,& \hspace{3mm}&\text{if}~u < \mu\\
 1&,& \hspace{3mm} &\text{if} ~ u \geq \mu
\end{array}
\right.
\end{eqnarray*}
where $\mu >0$ and $\eps >0$.\\\\
From the principle of conservation of volume ( developed in \cite{greenspan}), the dynamics of the tumor radius is governed by the following :
$$\frac{d}{d t} \Big(\frac{4 \pi R^3(t)}{3}\Big)= \int \int \int_{\Omega(t)} S(u)r^2 sin \theta~ d \theta~ d \phi~ dr-\int \int \int_{\Omega(t)}N(u) r^2 sin \theta~ d \theta~ d \phi~ dr,\hspace*{3mm} t>0$$
where $S(u)$ and $N(u)$ describe the proliferation and the mortality rates of tumor cells. They are given by
$$S(u)=\lam f(u),  \hspace*{10mm} N(u)=\eta >0.  $$
In general, the variational methods are the natural way to study the discontinuous elliptic problem, see for instance   \cite{badiale} and \cite{bonanno}, but with the classical variational methods, we can not characterize the variation of the free boundary. More precisely, let us introduce the notation
$$w^-(t)=\{u(r,t)<\mu \}, \hspace*{3mm} w^+(t)=\{u(r,t) \geq \mu\}$$
The boundary of $w^-$ denoted by $\partial w^-(t):=\Gamma(t)$ can be described by $r=r_0(t)$ for $t>0$ called the free boundary which is different from the tumor boundary $r=R(t)$. This type of problem is less studied in literature especially when $\Omega(t)$ is non spherical domain.\\\\
Denoting the initial tumor radius by $R_0$, the evolution of $R$ satisfies the following differential equation, for $t >0$
\begin{equation}\label{p3}
\left\{
\begin{array}{rcll}
\vspace*{2mm}R^2(t) \frac{d R(t)}{dt}&=&\int_{0}^{R(t)}  S(u)r^2dr-\int_{0}^{R(t)}N(u)r^2 dr\\
R(0)&=&R_0.\\
\end{array}
\right.
\end{equation}
When $\eps \in (0,1)$, in the region $w^+(t)$, nutrient concentration is enough to sustain tumor cells in normal proliferation. At some level $(u=\mu)$, the tumor cells are no longer able to continue at their normal pace and an inner region is developed namely $w^-(t)$ due to the nutrient deficiencies. Hence, we can resume the situation in the following problem
\begin{equation}\label{p4}
\left\{
\begin{array}{rcll}
\vspace*{3mm} \frac{1}{r^2} \frac{\partial}{\partial r} \big ( r^2 \frac{\partial u}{\partial r}\big)&=& \lam \big( \eps +(1-\eps ) H(u-\mu) \big) ~ \hspace{3mm}&~0< r <R(t),\hspace*{2mm} t >0\\
\vspace*{3mm} u(R(t),t)&=& u_{\infty},~\frac{\partial u}{\partial r}(0,t)=0 \hspace*{3mm}&~ \text{for} \hspace*{2mm} t >0\\
\vspace*{2mm}R^2(t) \frac{d R(t)}{dt}&=& \int_{0}^{R(t)} S(u) r^2 dr -\int_{0}^{R(t)} N(u) r^2 dr\\
R(0)&=&R_0.\\
\end{array}
\right.
\end{equation}
When $\eps=0$, we refer the reader to \cite{Byrne Medicine} for more details on the nonlinearity $f$.\\
When $\eps \in [1,+\infty)$, the study of problem (\ref{p4}) can be very interesting from a mathematical point of view. We refer to \cite{bensid-diaz1} and \cite{bensid-diaz2} for a similar case (in one dimension and without the variable $t$) for application to the Budyko climate models. \\\\
In this work, we are interested to study the existence, multiplicity of solutions of problem $(\ref{p4})$ and give the characterization of the free boundaries when the tumor region $\Omega(t)$ is non spherical domain.( More realistic).\\
The main difficulty is that the problem $(\ref{p4})$ has two different types of free boundaries. The outer boundary $\partial \Omega(t)$ ( moving boundary) and the inner boundary whose evolution is implicit generated by the discontinuous nonlinearity $H$( of obstacle-type). To the best of our knowledge, no investigation has been devoted to study problem  $(\ref{p4})$ when $\Omega(t)$ is a general domain. To give a positive answers, we use a perturbation method to reduces the study to solve a nonlinear integral equation and allow us to give a positive result to the solvability of the equation of free boundary under perturbation. Hence, having a solution of problem $(\ref{p4})$ in a spherical domain, we study the effect on the solutions under perturbation of the outer boundary. This approach was used by the second author for others free boundary problems ( one free boundary). See  \cite{Bensid4} and \cite{Bensid1}. \\\\
However, from the biological point of view, if the nutrient concentration becomes equal to another level (different from $\mu$), the tumor region $w^-(t)$ develops another region called necrotic region. To obtain this phase, we can focus the study on the following problem

\begin{equation}\label{pp}
	\left\{
	\begin{array}{rcll}
	\del u &=& \lam \displaystyle \sum_{i=1}^{2} \eps_i H(u-\mu_i)~ \hspace{3mm}&\text{in}~\Omega (t)\\
	u&=& \overline{u}_{\infty}~\hspace{3mm} &\text{on } ~ \partial \Omega(t)
	\end{array}
	\right.
\end{equation}
where $\eps_i \in (0,1)$.
\\There is also mathematical challenge in the study of the problem (\ref{pp}). The technical details and others results will be presented elsewhere.\\
We mention also that  when $f(u)=u,$ a perturbation analysis was used by Cui and Friedman in \cite{cui friedman} where the authors proved the existence of unique necrotic stationary solution and it is also asymptotically stable under small perturbation.  Friedman and Reitich \cite{FR} studied the problem (\ref{p2}), (\ref{p3}) and proved that there is a unique stationary solution for $u_{\infty} > \tilde{u}$. Others authors extended this result to general differentiable functions. We refer to \cite{cui}, \cite{FH}  and \cite{FR2}  and the reference therein for similar works.\\\\
 Recently, J. Wu and C. Wang \cite{wu-wang} and H. Song, W. Hu and Z. Wang \cite{wu3} studied problem (\ref{p4}) by considering the case $\eps \equiv 0$ in a spherical domain $\Omega(t)$ with Robin boundary condition (angiogenesis phenomenon). It is shown that there exists two thresholds, $\tilde{u}$ and $u^*$ giving the existence results. See \cite{wu3}. For others works on the existence of stationary solutions and the connection between the different phases, we refer to \cite{wu3}, \cite{Song1} and \cite{J. Wu}. \\\\
\hspace*{3mm} One of the main goal of this paper is to characterize the free boundary explicitly by obtaining the exact branches of radial solutions ( bifurcation diagram) and show the solvability of the problem in the neighborhood of the radial solutions under perturbation of the boundary value and a smooth boundary of the general domain $\Omega(t)$. Moreover, this perturbation procedure give a particular position of free boundary ( the inner boundary) which allows to generate a bifurcated solution of problem $(\ref{p4})$.\\ Obviously, by a solution we mean a couple $(u(r,t),R(t))$ satisfying the problem (\ref{p4}) and we assume the following assumption:\\\\
$(H_1)$ Assuming that the boundary $\partial \Omega(t)$ can be parameterised as $R(t)+\beta(\theta)$ for $t>0$ where $\beta \in C^2(S)$, $\theta \in S$ where $S$ is the unit sphere.\\\\
In this work, first, we will give the existence, multiplicity and diagram of bifurcation of solutions of problem $(\ref{p4})$ with some properties of their free boundaries when $\Omega(t)=B(0,R(t)).$\\
In this case, we show the existence of $r_0=r_0(t)$ such that the inner free boundary is a moving sphere given by $$u(r_0(t),\theta)=\mu,\quad t>0.$$
If we denote by $\chi_D$ the characteristic function of $D,$ and $w^-(t)$ the slowed growth region, then we prove that the following problem
\begin{equation}\label{chara}
\left\{
\begin{array}{rcll}
\del u &=& \lam(\eps +(1-\eps)) \chi_{\Omega(t) \setminus w^-(t)} ~ \hspace{3mm}&\text{in}~\Omega(t) \\
u&=& \overline{u}_{\infty}~\hspace{3mm} &\text{on } ~ \partial \Omega(t)
\end{array}
\right.
\end{equation}
has a unique solution. (See Proposition 4.1 below).\\
Then using the perturbation approach when $\Omega(t)$ satisfying $(H_1),$ we give the existence results of the two free boundaries of  problem $(1.1)$ ($D\equiv 1$). More precisely, if we convert the problem $(\ref{chara})$ into a nonlinear Hammerstein integral equation and prove the existence of function $b$ such that
$$u(r_0(t)+b(\theta),\theta)=\mu,\quad t>0,$$
then $u$ is solution of problem $(\ref{p1}).$\\
 We will prove also that in addition to radial solutions,  there exist a remarkable position of the free boundary corresponding to problem (\ref{p4}) when $\Omega(t)=B(0,R(t))$ for which a bifurcation phenomenon can occur. This situation is novelty in tumor growth models. \\\\
\hspace*{3mm} This paper is organized as follows : In section 2, we establish the existence and multiplicity of stationary solutions to the problem (\ref{p4}) when $\Omega=B(0,R)$. In section 3, we prove the global existence of transient solution. In section 4, we formulate  the perturbed problem and show that it is equivalent to a nonlinear integral equation. Section 5, concern the resolution of the integral equation by the local methods. Finally, an appendix is devoted to recall some useful results.

\section{Stationary solutions}
In this section, we study the stationary solutions to the problem (\ref{p4}), denotes by $(u_{\lambda,\mu},R)$. More precisely, we consider
\begin{equation}\label{ps}
\left\{\begin{array}{ll} \frac{1}{r^2} \frac{\partial}{\partial r} \big ( r^2 \frac{\partial u}{\partial r}\big)= \lam \big( \eps +(1-\eps ) H(u-\mu) \big),  & \quad \mbox{ }\ 0< r <R\\[0.3cm]
				u(R)= u_{\infty},\quad \quad \frac{\partial u}{\partial r}(0)=0        & \quad \mbox{ }\ \\[0.3cm]
				\frac{1}{R^2} \Big(\int_{0}^{R} S(u) r^2 dr -\int_{0}^{R}N(u) r^2 dr\Big) =0.       & \quad \mbox{ }\
			\end{array}
			\right.	
\end{equation}	

Recall that our nonnegative solutions must be strictly convex functions  such that $\displaystyle\min_{0<r<R}u_{\lambda,\mu}(r) =u_{\lambda,\mu}(0)$. The set $\Gamma=\{r_{\lam} \in (0,R),\hspace*{2mm} u_{ \lambda,\mu}(r_{\lam})=\mu\}$ is called the free boundary corresponding to the stationary solution.
\\In order to state the main result of this section, we define two crucial values of the
parameter $\lam$
$$ \lam_1:=\frac{6(u{_\infty}-\mu)}{R^2},\hspace*{3mm}
\lam_2:= \frac{27(u_{\infty}-\mu)(\eps-1)^2}{\eps^2 R^2 \big(\frac{4 \eps}{3}-\frac{3}{2}\big)}$$
The first theorem  concern the existence and multiplicity of stationary solutions of problem (\ref{ps}).
\begin{theo}\label{th1}
	\begin{enumerate}
	\item For $\eps \in (0,+\infty)$, we have the following:
	\begin{itemize}
		\item[$\bullet$] if $\lam < \lam_1$, then there exists a unique solution $u_{\lambda,\mu}$ without free boundary of (\ref{ps}). Moreover,
		$$u_{\lambda,\mu}(0)=\frac{-\lam}{6} R^2+u_{\infty}$$
		i.e. the line $(\lam, \gamma_s(\lambda)),\hspace*{3mm}\gamma_s(\lam):=\frac{-\lam}{6} R^2 +u_{\infty},\hspace*{3mm} \lam < \lam_1$
		\\defines a decreasing part of the bifurcation diagram.
	\end{itemize}
	\item For $\eps \in ( \frac{3}{2}, +\infty)$, we have the following:
	\begin{enumerate}
	\item[(i)] If $\lam=\lam_2$, then there exists a unique solution $u_{\lambda_2,\mu}^*$ of (\ref{ps}) giving rise to a free boundary given by
	$$r_{\lambda_2}=\frac{2 \eps-3}{3(\eps-1)} R$$
	\item[(ii)] If $\lam \in (\lam_2,\lam_1]$ then (\ref{ps}) has two distinct positive solutions $\overline{u}_{\lambda,\mu}$ and $\underline{u}_{\lambda,\mu}$ with theirs corresponding free boundaries $\overline{r}_\lam$ and $\underline{r}_\lam$ given explicitly by
	\begin{eqnarray*}
	\overline{r}_\lam&=&\frac{(\eps-\frac{3}{2})}{3(\eps-1)} R \Big[1+ 2cos\Big(\frac{1}{3}arccos\Big( 1+\frac{27(\eps-1)^2\big( \frac{1}{2} R^2-\frac{3(u_{\infty}-\mu)}{\lam}\big)}{2(\eps-\frac{3}{2})^3 R^2}\Big) \Big)\Big]\\
	\underline{r}_\lam&=&\frac{(\eps-\frac{3}{2})}{3(\eps-1)} R \Big[1+ 2cos\Big(\frac{1}{3}arccos\Big( 1+\frac{27(\eps-1)^2\big( \frac{1}{2} R^2-\frac{3(u_{\infty}-\mu)}{\lam}\big)}{2(\eps-\frac{3}{2})^3 R^2}\Big)+\frac{4 \pi}{3}\Big)\Big]
	\end{eqnarray*}
	Moreover;
	$$\overline{\gamma}_s(\lam)=-\frac{\lam \eps}{6} \overline{r}_\lam^2+\mu, \hspace*{6mm}\underline{\gamma}_s(\lam)=-\frac{\lam \eps}{6} \underline{r}_\lam^2+\mu$$
	\item [(iii)]If $\lam \in (\lam_1,+\infty)$, then there exists a unique positive solution $u_{\lambda,\mu}$ of (\ref{ps}) giving rise to a free boundary given by
	$$r_\lam=\frac{(\eps-\frac{3}{2})}{3(\eps-1)} R+\Big(\frac{-q+\sqrt{q^2+\frac{4}{27}p^3}}{2}\Big)^{\frac{1}{3}}+\Big(\frac{-q-\sqrt{q^2+\frac{4}{27}p^3}}{2}\Big)^{\frac{1}{3}}$$
	with
	$$p=\frac{-(\eps-\frac{3}{2})^2}{3(\eps-1)^2}R^2,\hspace*{3mm}
		q=\frac{-2(\eps-\frac{3}{2})^3 R^3}{27(\eps-1)^3}-\frac{ R \big( \frac{1}{2} R^2-\frac{3(u_{\infty}-\mu)}{\lam} \big)}{(\eps-1)}$$
	Moreover;
	$$\gamma_s(\lam)=-\frac{\lam \eps}{6} r_\lam^2+\mu$$
	\end{enumerate}
    \item For $\eps \in (0,\frac{9}{8}]$ and $\eps \neq 1$, if $\lam \geq \lam_1$ there exists one solution $u_{\lambda,\mu}$ of ( \ref{ps}) with a free boundary given by
    \begin{eqnarray*}
    r_\lam&=&\frac{(\eps-\frac{3}{2})}{3(\eps-1)} R \Big[1+ 2cos\Big(\frac{1}{3}arccos\Big( 1+\frac{27(\eps-1)^2\big( \frac{1}{2} R^2-\frac{3(u_{\infty}-\mu)}{\lam}\big)}{2(\eps-\frac{3}{2})^3 R^2}\Big) +\frac{4 \pi}{3}\Big)\Big]\text{when}~~ \eps<1\\
    r_\lam&=&\frac{(\eps-\frac{3}{2})}{3(\eps-1)} R \Big[1- 2cos\Big(\frac{1}{3}arccos\Big( 1+\frac{27(\eps-1)^2\big( \frac{1}{2} R^2-\frac{3(u_{\infty}-\mu)}{\lam}\big)}{2(\eps-\frac{3}{2})^3 R^2}\Big)\Big)\Big]\text{when}~~ \eps>1
    \end{eqnarray*}

    \item For $\eps \in (\frac{9}{8},\frac{3}{2})$, we have the following:
    \begin{enumerate}
    \item[(i)] If $\lam \in [\lam_1, \lam_2)$, then there one solution $u_{\lambda,\mu}$ of ( \ref{ps}) with a free boundary given by
    $$r_\lam=\frac{(\eps-\frac{3}{2})}{3(\eps-1)} R \Big[1- 2cos\Big(\frac{1}{3}arccos\Big( 1+\frac{27(\eps-1)^2\big( \frac{1}{2} R^2-\frac{3(u_{\infty}-\mu)}{\lam}\big)}{2(\eps-\frac{3}{2})^3 R^2}\Big)\Big)\Big]$$
    \item[(ii)] If $\lam \in [\lam_2,+\infty)$, then there one solution $u_{\lambda,\mu}$ of ( \ref{ps}) with a free boundary given by
    $$r_\lam=\frac{(\eps-\frac{3}{2})}{3(\eps-1)} R+\Big(\frac{-q+\sqrt{q^2+\frac{4}{27}p^3}}{2}\Big)^{\frac{1}{3}}+\Big(\frac{-q-\sqrt{q^2+\frac{4}{27}p^3}}{2}\Big)^{\frac{1}{3}}$$
    \end{enumerate}
  \item If $\eps=\frac{3}{2}$, for $\lam \in (\lam_1,+\infty)$, then there exists a unique positive solution $u_{\lambda,\mu}$ of (\ref{ps}) giving rise to a free boundary given by
  $$r_\lam=\Big[R\Big(R^2-\frac{6(u_{\infty}-\mu)}{\lam}\Big)\Big]^{\frac{1}{3}}$$
	\end{enumerate}
\end{theo}
 \begin{figure}[h]
	\includegraphics[width=8 cm , height=5 cm]{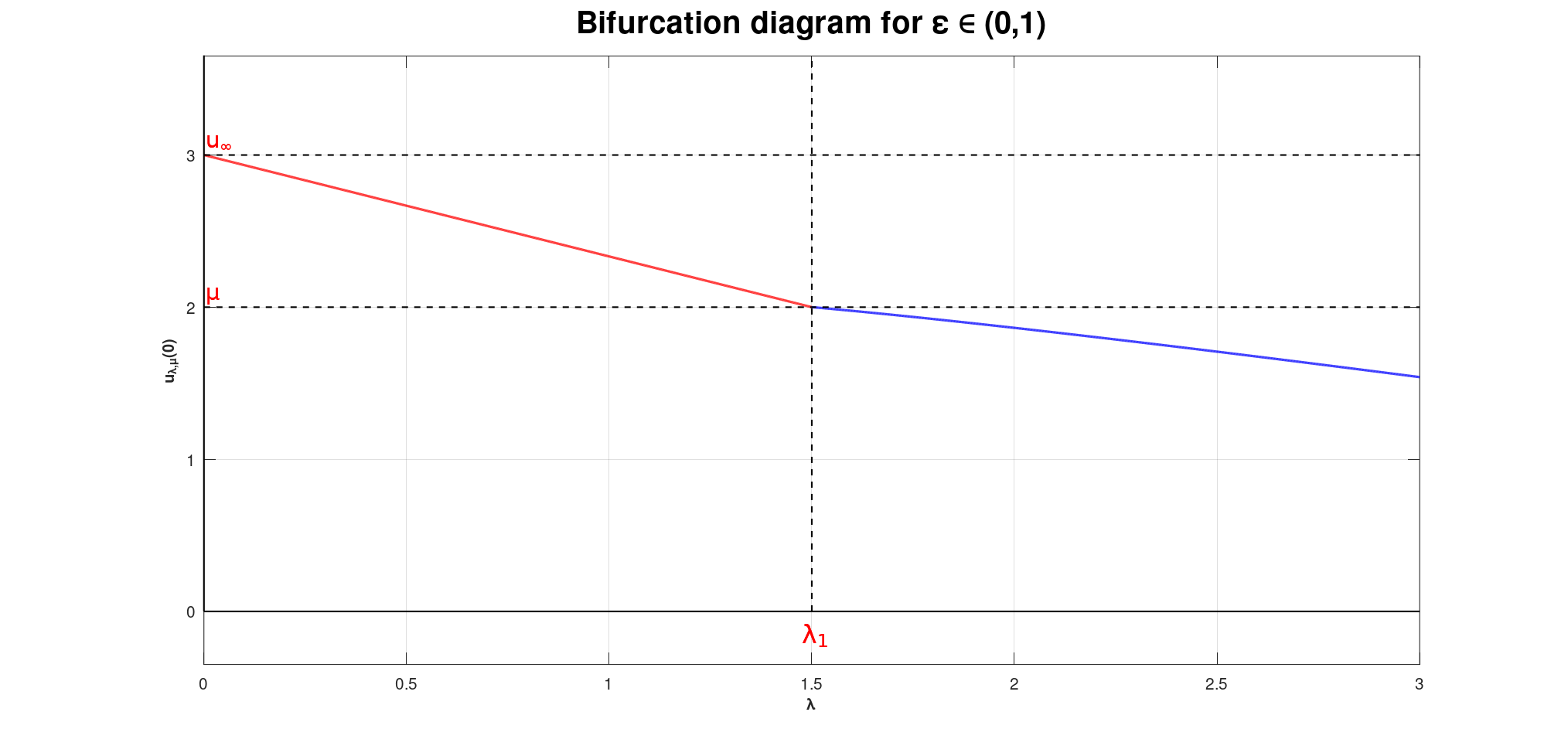}
	\includegraphics[width=9 cm , height=5 cm]{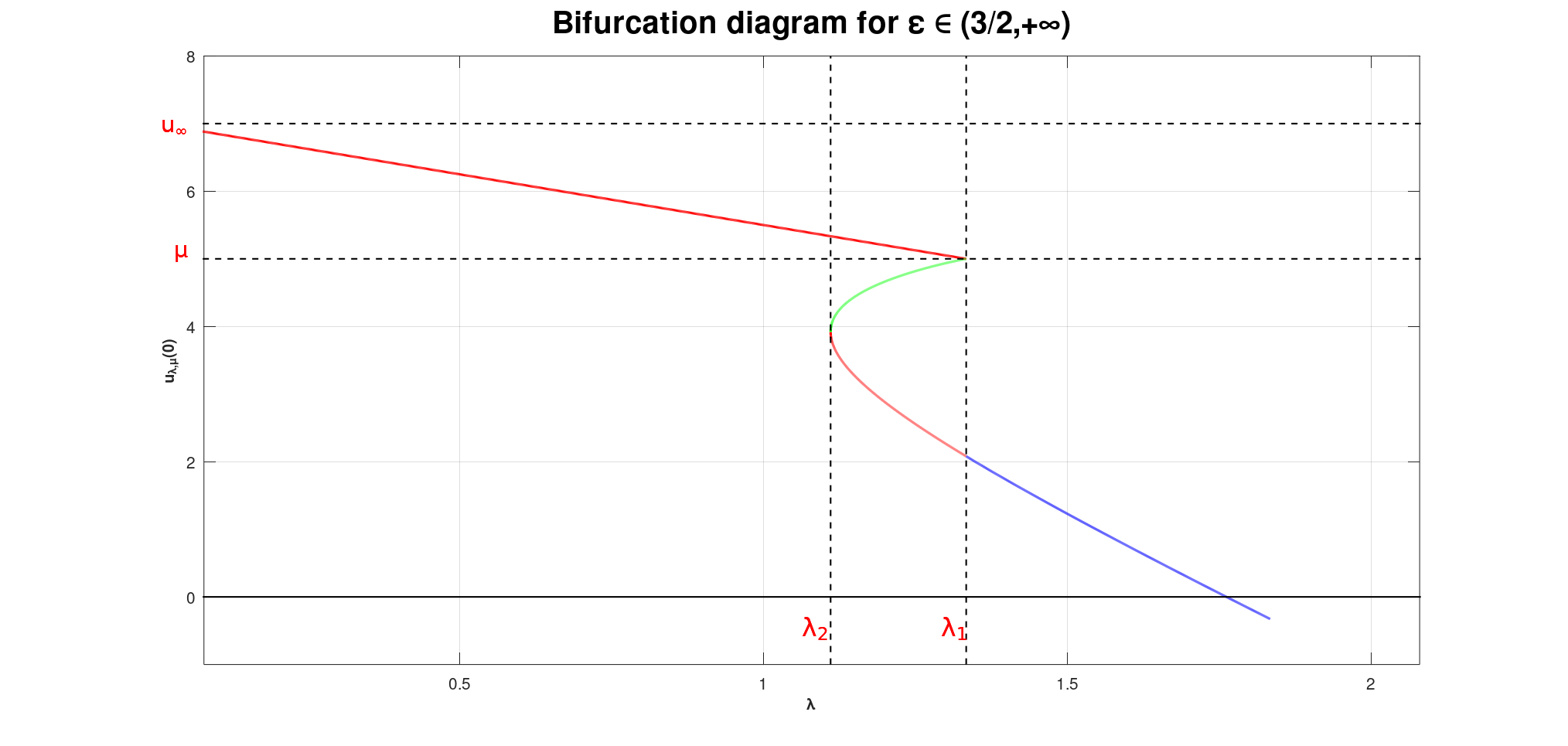}
\begin{center}
	\includegraphics[width=9 cm , height=5 cm]{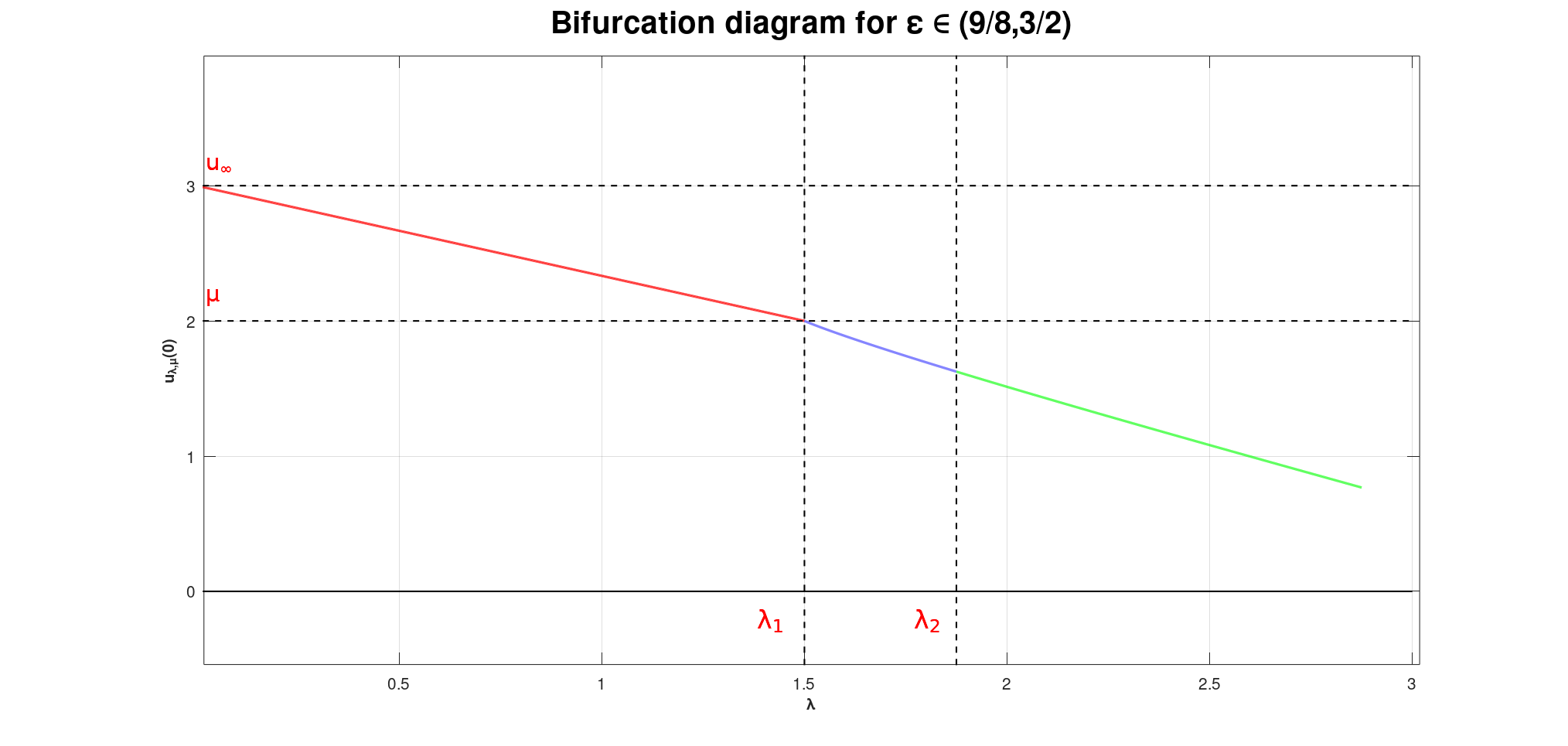}
\end{center}
\end{figure}
\begin{rem}
\begin{itemize}
\item The previous $\lam$ diagrams of bifurcation $(\lam,u_{\lambda,\mu}(0))$ give the exact numbers of solution of problem (\ref{ps}) which depend on the position of parameter $\lam=\lam(\eps,\mu,R)$. For example, when $\eps > 3/2$ and $\lam_2 <\lam <\lam_1$, there exists two solutions with free boundaries and one (in red) without free boundary. The explicit value of free boundaries stated in Theorem \ref{th1} and the different branches of bifurcation curves can be simplify the study of stability of our solutions. The technical details will be presented elsewhere.
\item When $\eps=1,$ problem $(1.4)$ becomes a classical elliptic problem without free boundary.
\end{itemize}
\end{rem}
\prrr 21
\\First, we consider the case $u_{\lambda,\mu}(0)>\mu$ corresponding to the case without free boundary. We have the problem
\begin{equation*}
\hspace*{0.5cm}
\left\{
\begin{array}{rcll}
\vspace*{5mm} \frac{1}{r^2} \frac{\partial}{\partial r} \big ( r^2 \frac{\partial u}{\partial r}\big)&=& \lam~ \hspace{3mm}&~0< r <R\\
u(R)=u_{\infty} &,&\hspace*{3mm} u'(0)=0
\end{array}
\right.
\end{equation*}
The solution is given by
$$u_{\lambda,\mu}(r)=\frac{\lam}{6}\big(r^2-R^2\big) +u_{\infty},\hspace*{2mm} r \in [0,R]$$
Since $$\displaystyle \min_{0<r<R} u_{\lambda,\mu}(r)=u_{\lambda,\mu}(0)=-\frac{\lam}{6} R^2+u_{\infty} > \mu~~ \text{if and only if }~~\lam < \frac{6(u_{\infty}-\mu)}{R^2} $$
Hence, this prove (1) of Theorem \ref{th1}.\\
\hspace*{3mm} Now, we consider the case $u_{\lambda,\mu}(0) \leq \mu$ with one free boundary, we look for the free boundary in the form $\{(r_\lam,\tet), \tet \in S\}$ for some $r_\lam \in (0,R)$ with $u(r_\lam)=\mu$, so that the differential equation
$$\frac{1}{r^2} \frac{\partial}{\partial r} \big ( r^2 \frac{\partial u}{\partial r}\big)= \lam \eps ~ \hspace{3mm}~0< r <r_\lam$$
$$\frac{1}{r^2} \frac{\partial}{\partial r} \big ( r^2 \frac{\partial u}{\partial r}\big)= \lam  ~ \hspace{3mm}~r_\lam< r <R$$
 $$u(R)=u_{\infty} ,\hspace*{3mm} u'(0)=0$$
with the following transmission condition on the free boundary
$$\frac{\partial u}{\partial r} \big(r_\lam^-\big)=\frac{\partial u}{\partial r} \big(r_\lam^+\big)$$
where $\frac{\partial u}{\partial r} \big(r_\lam^-\big)$ denotes the left derivative of $u$ and $\frac{\partial u}{\partial r} \big(r_\lam^+\big)$ denotes the right derivative at the value $r = r_\lam$.
\\In the region $(0,r_\lam)$, we have
\begin{equation*}
\hspace*{0.5cm}
\left\{
\begin{array}{rcll}
 \vspace*{5mm} \frac{1}{r^2} \frac{\partial}{\partial r} \big ( r^2 \frac{\partial u}{\partial r}\big)&=& \lam \eps~ \hspace{3mm}&~ r \in (0,r_\lam)\\
u(r_\lam)=\mu &,&\hspace*{3mm} u'(0)=0
\end{array}
\right.
\end{equation*}
and then,
$$u_{\lambda,\mu}(r)=\frac{\lam \eps}{6} \big( r^2-r_\lam^2 \big) +\mu,\hspace*{3mm} 0<r<r_\lam$$
In the region $(r_\lam,R)$, we have the following
\begin{equation*}
\hspace*{0.5cm}
\left\{
\begin{array}{rcll}
\vspace*{5mm} \frac{1}{r^2} \frac{\partial}{\partial r} \big ( r^2 \frac{\partial u}{\partial r}\big)&=& \lam ~ \hspace{3mm}&~ r \in (r_\lam,R)\\
u(r_\lam)=\mu &,&\hspace*{3mm} u(R)= u_{\infty}
\end{array}
\right.
\end{equation*}
and thus,
$$u_{\lambda,\mu}(r)=u_{\infty}+\frac{\lam}{6} \big( r^2-R^2 \big) +\frac{(r-R)r_\lam}{(r_\lam-R)r} \big(\mu - u_{\infty}-\frac{\lam}{6} \big( r_\lam^2-R^2 \big) \big),\hspace*{3mm} r_\lam<r<R$$
Using the transmission condition, implies that necessarily
\begin{equation}\label{11}
	\lam=\frac{3(u_{\infty}-\mu)}{(\eps-1)r_\lam^2 \frac{(R-r_\lam)}{R}+\frac{1}{2}(R^2-r_\lam^2)}
\end{equation}
In order to study this condition, we introduce the auxiliary function
\begin{equation}
	g(r)=\frac{3(u_{\infty}-\mu)}{(\eps-1)r^2 \frac{(R-r)}{R}+\frac{1}{2}(R^2-r^2)}~~\text{for}~~r \in (0,R)
\end{equation}
With
$$g'(r)=\frac{-3(u_{\infty}-\mu)r \big(\frac{\eps-1}{R}(2R-3r)-1\big)}{\big(r^2(\eps-1)\frac{R-r}{R}+\frac{1}{2}(R^2-r^2)\big)^2}$$
To study the behavior of function $g$, we distinguish two cases :
\begin{enumerate}
\item For $\eps \in (3/2,+\infty)$,
 the function has a global minimum at $r_{\lam_2}$ given by:
$$r_{\lam_2}=\frac{2\eps-3}{3(\eps-1)}R,\hspace*{3mm} g(r_{\lam_2})=\frac{27(u_{\infty}-\mu)(\eps-1)^2}{\eps^2 R^2 \big(\frac{4 \eps}{3}-\frac{3}{2}\big)}$$	
Therefore, when $\lam =\lam_2:=g(r_{\lam_2})$ the (\ref{ps}) has a unique solution, and the free boundary is given by $r_{\lam_2}$.
\\  When $\lam_2 <\lam\leq \lam_1:=g(0)$, it follows that  (\ref{11}) has exactly two roots between $(0,R)$, and when $\lam >\lam_1$, the equation (\ref{11}) has one roots.
\item For $\eps \in (0,\frac{3}{2}]$,
the function $g$ is increasing and has a minimum at $0$ given by $\frac{6(u_{\infty}-\mu)}{R^2}$. Thus, for $\lam \geq \lam_1$, the equation (\ref{11}) has one roots.
\end{enumerate}
\begin{figure}[h]
	\includegraphics[width=8 cm , height=5 cm]{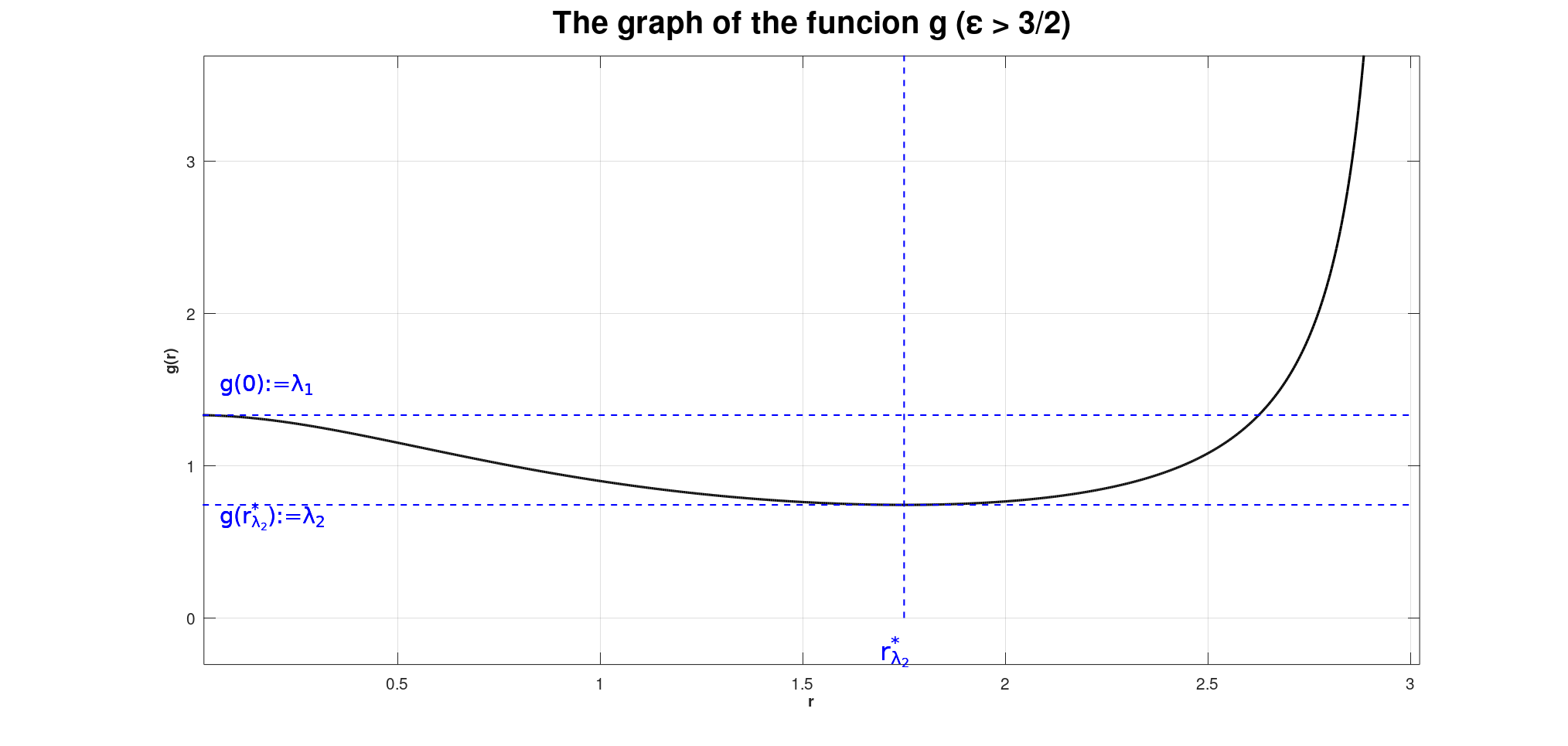}
	\includegraphics[width=8 cm , height=5 cm]{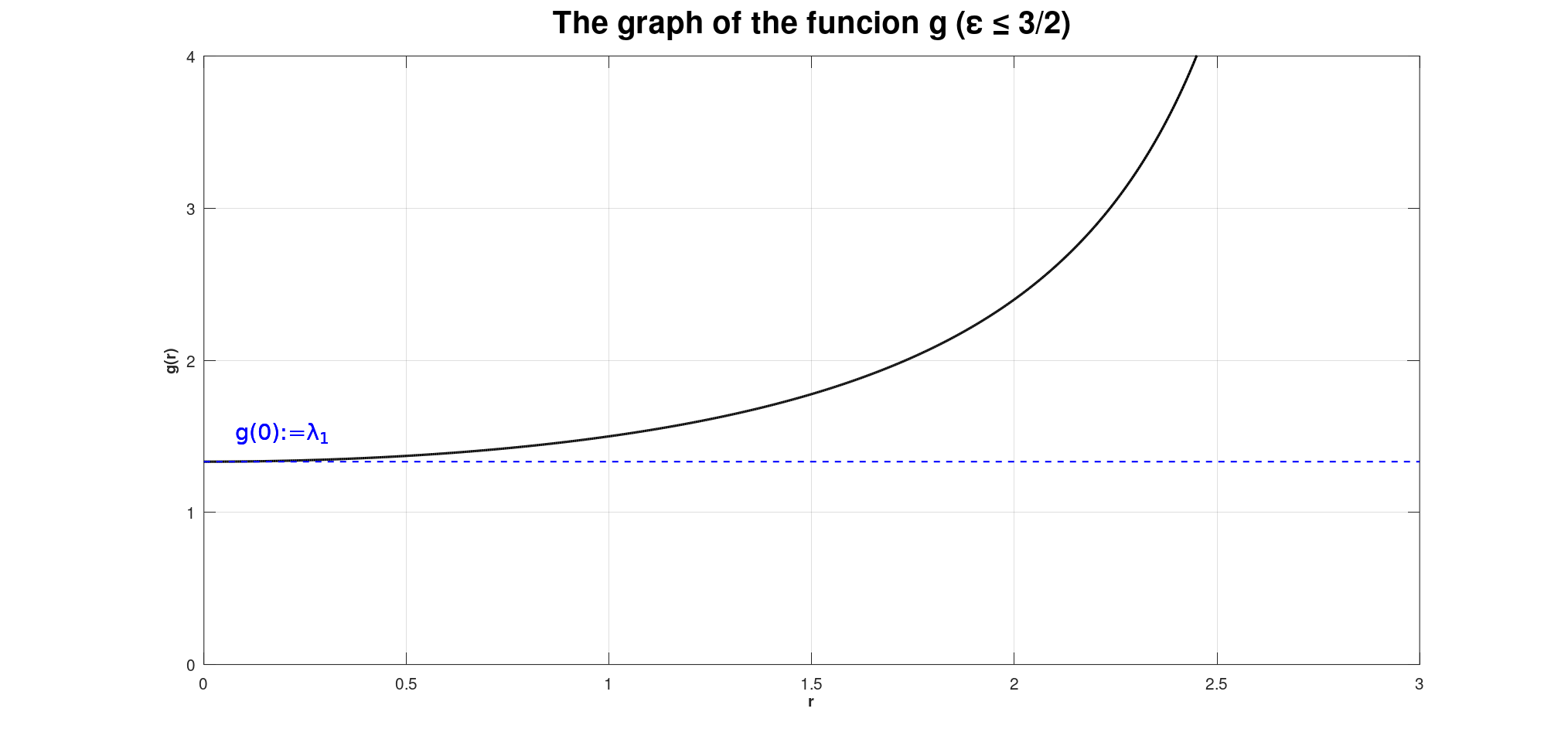}
\end{figure}
We obtain explicitly the free boundaries by solving the follow algebraic equation
\begin{equation}\label{12}
\frac{-(\eps-1)}{R} r^3+(\eps-\frac{3}{2}) r^2+ \frac{1}{2} R^2-\frac{3(u_{\infty}-\mu)}{\lam}=0,\hspace*{3mm} 0<r<R
\end{equation}
Considering the discriminant of this equation given by $\Delta = q^2+\frac{4}{27} p^3$,
with
$$  p=\frac{-(\eps-\frac{3}{2})^2}{3(\eps-1)^2}R^2, \hspace*{3mm}
	q=\frac{-2(\eps-\frac{3}{2})^3 R^3}{27(\eps-1)^3}-\frac{ R \big( \frac{1}{2} R^2-\frac{3(u_{\infty}-\mu)}{\lam} \big)}{(\eps-1)}.$$
Hence,
$$\Delta =\Big(\frac{R^2(\frac{1}{2}R^2-\frac{3(u_{\infty}-\mu)}{\lam})}{2(\eps-1)^2}\Big)\Big(\frac{1}{2}R^2+\frac{4(\eps-\frac{3}{2})^3}{27(\eps-1)^2}R^2-\frac{3(u_{\infty}-\mu)}{\lam}\Big)$$
(see Appendix A ).
\\ \\
If $\eps >3/2$, we have
\begin{itemize}
\item[$\bullet$] When $\lam_2 <\lam\leq \lam_1$, $\Delta <0$ and the two roots are given by
\begin{eqnarray*}
	\overline{r}_\lam&=&\frac{(\eps-\frac{3}{2})}{3(\eps-1)} R \Big[1+ 2cos\Big(\frac{1}{3}arccos\Big( \frac{-q}{2}\sqrt{\frac{27}{-p^3}} \Big) \Big)\Big]\\
	&=&\frac{(\eps-\frac{3}{2})}{3(\eps-1)} R \Big[1+ 2cos\Big(\frac{1}{3}arccos\Big( 1+\frac{27(\eps-1)^2\big( \frac{1}{2} R^2-\frac{3(u_{\infty}-\mu)}{\lam}\big)}{2(\eps-\frac{3}{2})^3 R^2}\Big) \Big)\Big]
\end{eqnarray*}
and
\begin{eqnarray*}
	\underline{r}_\lam&=&\frac{(\eps-\frac{3}{2})}{3(\eps-1)} R \Big[1+ 2cos\Big(\frac{1}{3}arccos\Big( \frac{-q}{2}\sqrt{\frac{27}{-p^3}} \Big)+\frac{4 \pi}{3} \Big)\Big]\\
	&=&\frac{(\eps-\frac{3}{2})}{3(\eps-1)} R \Big[1+ 2cos\Big(\frac{1}{3}arccos\Big( 1+\frac{27(\eps-1)^2\big( \frac{1}{2} R^2-\frac{3(u_{\infty}-\mu)}{\lam}\big)}{2(\eps-\frac{3}{2})^3 R^2}\Big)+\frac{4 \pi}{3}\Big)\Big]
\end{eqnarray*}
\item[$\bullet$] When $\lam >\lam_1$, $\Delta>0$ and we obtain the following roots
$$r_\lam=\frac{(\eps-\frac{3}{2})}{3(\eps-1)} R+\Big(\frac{-q+\sqrt{\Delta}}{2}\Big)^{\frac{1}{3}}+\Big(\frac{-q-\sqrt{\Delta}}{2}\Big)^{\frac{1}{3}}.$$
\end{itemize}
To get ((i), (ii) and (iii)) of theorem (\ref{th1}), it suffices to denote by $u_{\lam_2,\mu}^*$, $\overline{u}_{\lam,\mu}$, $\underline{u}_{\lam,\mu}$ and $u_{\lam,\mu}$ the solutions associated to the free boundaries given by $r_{\lam_2}$, $\overline{r}_\lam$, $\underline{r}_\lam$ and $r_\lam$ respectively.
\\ \\
If $\eps \in (0,\frac{3}{2}]$, we distinguish three cases
\begin{itemize}
\item[*] For $\eps \in (0,\frac{9}{8}]$, when $\lam >\lam_1$, its roots are the following
$$r_\lam=\frac{(\eps-\frac{3}{2})}{3(\eps-1)} R \Big[1+ 2cos\Big(\frac{1}{3}arccos\Big( 1+\frac{27(\eps-1)^2\big( \frac{1}{2} R^2-\frac{3(u_{\infty}-\mu)}{\lam}\big)}{2(\eps-\frac{3}{2})^3 R^2}\Big)+\frac{4 \pi}{3}\Big)\Big]\text{when}~~ \eps<1$$
 $$r_\lam=\frac{(\eps-\frac{3}{2})}{3(\eps-1)} R \Big[1- 2cos\Big(\frac{1}{3}arccos\Big( 1+\frac{27(\eps-1)^2\big( \frac{1}{2} R^2-\frac{3(u_{\infty}-\mu)}{\lam}\big)}{2(\eps-\frac{3}{2})^3 R^2}\Big)\Big)\Big]\text{when}~~ \eps>1$$
\item[*] For $\eps \in (\frac{9}{8},\frac{3}{2})$, we have
\\$\bullet$ When $\lam \in [\lam_1,\lam_2)$, the roots is given by
 $$r_\lam=\frac{(\eps-\frac{3}{2})}{3(\eps-1)} R \Big[1- 2cos\Big(\frac{1}{3}arccos\Big( 1+\frac{27(\eps-1)^2\big( \frac{1}{2} R^2-\frac{3(u_{\infty}-\mu)}{\lam}\big)}{2(\eps-\frac{3}{2})^3 R^2}\Big)\Big)\Big]$$
 $\bullet$ When $\lam \geq \lam_2$, we obtain the following roots
 $$r_\lam=\frac{(\eps-\frac{3}{2})}{3(\eps-1)} R+\Big(\frac{-q+\sqrt{\Delta}}{2}\Big)^{\frac{1}{3}}+\Big(\frac{-q-\sqrt{\Delta}}{2}\Big)^{\frac{1}{3}}$$
 \item[*] If $\eps= \frac{3}{2}$, for $\lam \in [\lam_1,+\infty)$ we have $\Delta=(q)^2 >0$ and the roots is given by
  $$r_\lam=\Big(\frac{-q+\sqrt{\Delta}}{2}\Big)^{\frac{1}{3}}+\Big(\frac{-q-\sqrt{\Delta}}{2}\Big)^{\frac{1}{3}}=\Big[R\Big(R^2-\frac{6(u_{\infty}-\mu)}{\lam}\Big)\Big]^{\frac{1}{3}}$$
\end{itemize}

Thus, the proof of Theorem \ref{th1} ends.\cqd
Now from theorem \ref{th1}, we see that we have two phases, $u_{\lam,\mu}(0) > \mu$ and $u_{\lam,\mu}(0) \leq \mu$. This leads to determine the radius at which the tumor moves from phase 1 ($u_{\lam,\mu}(0) > \mu$) into phase 2 ($u_{\lam,\mu}(0) \leq \mu$) when $\eps \in (0,1)$.\\
\\So, the tumor remains in phase 1 until $u_{\lam,\mu}(0) = \mu$ ( i.e
$\frac{-\lam}{6} R^2+u_{\infty}=\mu$),
giving the outer radius by $$R^*=\sqrt{\frac{6(u_{\infty}-\mu)}{\lam} }.$$
When $R < R^*$, the tumor region is in phase 1 (without free boundary) and when $R\geq R^*$, the tumor develops a new region (with free boundary).
\\In other part, to determine the existence of $R$ in stationary case, we need to find $R$ satisfying the following equation
$$\int_{0}^{R} S(u) r^2 dr -\int_{0}^{R}N(u) r^2 dr=0$$
where $S(u)=\lam f(u)$ and $N(u)=\eta$. Let
\begin{eqnarray*}
L(R)&:=&\lam \int_{0}^{R} f(u) r^2 dr -\int_{0}^{R}N(u) r^2 dr=\int_{0}^{r_\lam} \lam \eps r^2 dr +\int_{r_\lam}^{R} \lam r^2 dr -\int_{0}^{R} \eta r^2 dr\\
&=&\frac{\lam (\eps-1)}{3} r_\lam^3 + \frac{\lam-\eta}{3} R^3, \hspace*{3mm}\text{for} ~~R >0,
\end{eqnarray*}
where $r_\lam$ is given in theorem \ref{th1} when $\eps \in (0,1)$. Hence, $L(R)=0$, for $R>0$ implies that
$$\frac{\lam (\eps-1)}{3} \Big(\frac{r_\lam}{R} \Big)^3 +\frac{\lam-\eta}{3}=0.$$
Assuming that $\lam > \eta $, then we have :
\begin{itemize}
\item For $R <R^*$, the tumor remains in phase 1 and $r_\lam=0$. So
$$\lim_{R \to 0} \frac{L(R)}{R^3}=\frac{\lam-\eta}{3} >0.$$
\item For $R \geq R^*$, the tumor is in phase 2. The behavior of function $g$ ( see (2.9)) give that $\lam=\lam(R)$ is non increasing with respect to $R,$ then $\displaystyle \lim_{R \to +\infty}\lam(R)=0.$ Hence, 
$$\lim_{R \to +\infty} \frac{L(R)}{R^3}=-\frac{\eta}{3} <0.$$
\end{itemize}
A standard argument shows  the existence of $R:=R_s$ such that $L(R_s)=0$.
\section{Transient solutions}
In this section, we study the existence of global solutions and asymptotic behavior of solution of free boundary problem (\ref{p4}).
The first result is the following
\begin{theo}\label{th31}
For any $R_0 >0$ and for $0< \eps <1$, the problem (\ref{p4}) has a unique global solution $(u(r,t),R(t))$ for $t >0$.
\end{theo}
\begin{rem}
	The result stated in theorem \ref{th31} remains true for all $\eps >0$ excepted for $\eps >3/2$ and $\lam \in ]\lam_2,\lam_1]$. In this case, the uniqueness fails (see $2,ii)$ from theorem \ref{th1}). Moreover, the case $ \eps >1$ has no biological meaning in tumor growth models.
\end{rem}
For $t >0$,
\begin{eqnarray*}
\frac{d R(t)}{d t}&=& \frac{1}{R(t)^2} \Big( \int_{0}^{r_0(t)} \lam \eps r^2 dr +\int_{r_0(t)}^{R(t)} \lam r^2 dr -\int_{0}^{R(t)} \eta r^2 dr \Big)=\frac{\lam (\eps-1)}{3} \frac{r_0(t)^3}{R(t)^2} + \frac{\lam-\eta}{3} R(t)\\
&=& R(t) \Big[\frac{\lam (\eps-1)}{3} \Big(\frac{r_0(t)}{R(t)} \Big)^3 + \frac{\lam-\eta}{3} \Big]
\end{eqnarray*}
Then, we get a differential equation which governs the motion radius $R$ depending on the radius $r_0(t)=r_{\lam}(R(t))$, where $r_{\lam}$ is given by the theorem \ref{th1} when $\eps \in (0,1)$.
\\We set
\begin{equation}\label{fonctionh}
H(R)= \frac{\lam (\eps-1)}{3} \Big(\frac{r_0(t)}{R(t)} \Big)^3 +\frac{\lam-\eta}{3}
\end{equation}
where
$$\frac{r_0(t)}{R(t)}=\frac{(\eps-\frac{3}{2})}{3(\eps-1)} \Big[1+ 2cos\Big(\frac{1}{3}arccos\Big( 1+\frac{27(\eps-1)^2\big( \frac{1}{2} R(t)^2-\frac{3(u_{\infty}-\mu)}{\lam}\big)}{2(\eps-\frac{3}{2})^3 R(t)^2}\Big) +\frac{4 \pi}{3}\Big)\Big]$$
The proof of theorem \ref{th31} is based on the follow lemma
\begin{lem}\label{lemh}
For $\eps \in (0,1)$, the function $H$ is decreasing and satisfies
\begin{equation}\label{h}
	\frac{\lam \eps -\eta}{R} \leq H(R) \leq \frac{\lam-\eta}{3}
\end{equation}
\end{lem}
\pr
An easy calculation shows that;
$$\frac{\partial }{\partial R}\Big(\frac{r_0}{R} \Big) =\frac{3(\eps-1)}{(\eps-\frac{3}{2})^2}~ \frac{u_{\infty}-\mu}{\lambda R^3(t) \sqrt{1-A^2(R)}}~ sin\Big(\frac{1}{3} arccos\big(A(R)\big)+\frac{4 \pi}{3}\Big)~~>0$$
where $$A(R)=1+\frac{27(\eps-1)^2\big( \frac{1}{2} R(t)^2-\frac{3(u_{\infty}-\mu)}{\lam}\big)}{2(\eps-\frac{3}{2})^3 R(t)^2}$$
then,
$$H'(R)=\lam (\eps-1) ~\Big(\frac{r_0}{R} \Big)^2 ~\frac{\partial }{\partial R}\Big(\frac{r_0}{R} \Big)~~<0 \hspace*{3mm}\text{for}~\eps \in (0,1)$$
Hence, $H$ is strictly decreasing.
\\For $\eps \in (0,1)$ and $0 \leq r_0 \leq R$, we have
$$\frac{\lam (\eps-1)}{3} \leq \frac{\lam (\eps-1)}{3} \Big(\frac{r_0}{R} \Big) ^3 \leq 0$$
then, we obtain explicitly (\ref{h}).This proves lemma \ref{lemh}. \cqd
\prrr{3.1} Given $R(t)>0$ and under the condition $0<\eps <1$, we know that the problem (\ref{p4}) admits a unique solution $u(r,R(t))$. So, we can determine $R(t)$ by solving the following initial problem
\begin{equation}
\hspace*{0.5cm}
\left\{
\begin{array}{rcll}
R'(t) &=& R(t) H(R(t))~ \hspace{3mm}&t >0\\
R(0)&=& R_0
\end{array}
\right.
\end{equation}
where $H$ is given by (\ref{fonctionh}).
\\ The result of lemma \ref{lemh} imply that for $R_0 >0$, problem (\ref{p4}) has a unique global solution satisfying :

\begin{eqnarray}
R_0 e^{\big(\frac{\lam \eps-\eta}{3}\big)t} \leq R(t) \leq R_0 e^{\big(\frac{\lam-\eta}{3}\big)t}&,& \hspace*{3mm}t>0.\label{r3}
\end{eqnarray}
This conclude the proof of Theorem \ref{th31}.\cqd
Next, the result concerning the asymptotic behavior of transient solution, we have
\begin{theo}\label{th32}
For any initial value $R_0>0$ and for $0<\eps <1$, we have
\begin{enumerate}
\item If $\eta > \lam$, then $\displaystyle\lim_{t\to +\infty} R(t)=0$
\item If $\eta \leq \lam $, then
$$\lim_{t \to +\infty} R(t)=R_s \hspace*{4mm} \lim_{t \to +\infty} u(r,t)=u_{\lam,\mu}$$
where $(u_{\lam,\mu},R_s)$ is the stationary solution of the problem (\ref{p4}).
\end{enumerate}
\end{theo}
\prrr{3.2} If $\lam > \eta$, then from (\ref{r3}),
we conclude immediately that \\ $\displaystyle\lim_{t\to +\infty} R(t)=0$.
\\\\If $\lam \leq \eta$, we have from lemma \ref{lemh} $H$ is decreasing function and the equation $H(R)=0$  has a unique positive constant solution $R_s$  ( see the last equation of problem (\ref{ps})), then a classical result of differential equation theory give that
$$\lim_{t \to +\infty} R(t)=R_s \hspace*{3mm}\text{for any initial value} ~R_0 >0.$$
So, $\displaystyle\lim_{t \to +\infty} u(r,t)=u_{\lam,\mu}.$
\cqd
\section{The perturbed problem}
\hspace*{3mm} In this section, we are concerned with problem (\ref{p1}) when $\Omega(t)$ verifies $(H_1)$. Recall  that the boundary $\partial \Omega(t)$  can be parameterized as $R(t)+\bet(\tet)$, where $\bet\in C^2(S)$ and $S$ is the unit sphere.
\\\\We denote by $\Omega_{\bet}(t)$ the admissible perturbation of the ball $B(0,R(t))$. So,
$$\Omega(t):=\Omega_{\bet}(t)=B(0,R(t)+\bet(\tet)),\hspace*{2mm} \Omega_{0}(t)=B(0,R(t)).$$ Let $r_{0}$ denote one of the values of the free boundary considering in Theorem \ref{th1}. In view of Theorem \ref{th32} we can consider only the stationary case. Otherwise, we fix any $t=T$ and we look for the free boundary in the form $r_0+b(\theta),$  where $b $ is the perturbation caused by $\beta.$ Hence, we consider the following problem
\begin{equation}\label{sbeta}
\left\{
\begin{array}{rcll}
\del u &=& \lam ( \eps +(1-\eps ) H(u-\mu))~ \hspace{3mm}&\text{in}~\Omega_{\bet} \\
u&=& \overline{u}_{\infty}~\hspace{3mm} &\text{on } ~ \partial \Omega_{\bet}
\end{array}
\right.
\end{equation}
where $\Omega_{\bet}$ is the admissible perturbation of the ball $B(0,R)$; i.e $\Omega_{\bet} =B(0,R+\bet(\theta))$ and $\overline{u}_{\infty}$ is closed to $u_{\infty}$.
\\We define the set of admissible surfaces in $\Omega_{\bet}$ by
$$S_{\beta}=\{f \in C(S); (f(\theta),\theta) \in \Omega_{\bet} \hspace*{2mm}for\hspace*{2mm} \theta \in S\}$$
For a function $\psi \in S_{\beta}$, we consider the set :
$$\Omega_{\bet,\psi}=\{(r,\theta) \in \Omega_{\bet},\hspace*{3mm} r<\psi(\theta)\}$$
We denote by $\chi_{\Omega_{\bet} \setminus \Omega_{\bet,\psi}}$ the characteristic function of $\Omega_{\bet} \setminus \Omega_{\bet,\psi}$ then, we have the following result.
\begin{prop}
Assume that:
\begin{eqnarray*}
(H_2)\hspace*{3mm}\lam &\geq& \lam_2, \hspace*{3mm}for\hspace*{2mm} \eps \in (3/2,+\infty)\\
(H_3)\hspace*{3mm} \lam &\geq& \lam_1, \hspace*{3mm}for\hspace*{2mm} \eps \in (0,3/2]
\end{eqnarray*}
Then the problem
\begin{equation}\label{betaxis}
\left\{
\begin{array}{rcll}
\del u &=& \lam(\eps +(1-\eps) \chi_{\Omega_{\bet} \setminus \Omega_{\bet,\psi}}) ~ \hspace{3mm}&\text{in}~\Omega_{\bet}\\
u&=& \overline{u}_{\infty}~\hspace{3mm} &\text{on } ~ \partial \Omega_{\bet}
\end{array}
\right.
\end{equation}
has a unique solution $u_{\lambda,\beta} \in C^{1,\alpha}(\overline{\Omega_{\bet}}, \RR)$ with $\alpha=1-\frac{3}{p}$. Moreover, if $u_{\lambda,\beta}(\psi(\theta),\theta)=\mu$ and $\overline{u}_{\infty} > \mu$ then $u_{\lambda,\beta}$ is solution of (\ref{sbeta}).
\end{prop}
\prr{3.1}
\\We see that $\lam(\eps +(1-\eps) \chi_{\Omega_{\bet} \setminus \Omega_{\bet,\psi}}) \in L^p(\Omega_{\bet})$, $p >1$. From \cite{elliptic}, there exists a unique solution of (\ref{betaxis}) in $W^{2,p}(\Omega_{\bet})$. For $p >3$, $W^{2,p}(\Omega_{\bet}) \subset C^{1,\alpha}(\overline{\Omega_{\bet}}, \RR)$ with $\alpha=1-\frac{3}{p}$.
\\The solution $u_{\lambda,\beta}$ satisfies;
\begin{eqnarray*}
\del u_{\lambda,\beta}&=& \lam \eps , \hspace*{3mm} \Omega_{\bet,\psi}\\
\del u_{\lambda,\beta}&=& \lam  , \hspace*{3mm} \Omega_{\bet} \setminus \Omega_{\bet,\psi}\\
u_{\lambda,\beta}&=& \overline{u}_{\infty} , \hspace*{3mm} \partial \Omega_{\bet}
\end{eqnarray*}
If we prove the existence of a function $\psi$ such that $u_{\lambda,\beta}(\psi(\theta),\theta)=\mu$, then $u_{\lambda,\beta}$ will be a solution of
\begin{eqnarray*}
	\del u_{\lambda,\beta}&=& \lam \eps , \hspace*{3mm} \Omega_{\bet,\psi}\\
	u_{\lambda,\beta}&=& \mu ,\hspace*{3mm} \partial \Omega_{\bet,\psi}\\
	\del u_{\lambda,\beta}&=& \lam  , \hspace*{3mm} \Omega_{\bet} \setminus \Omega_{\bet,\psi}\\
	u_{\lambda,\beta}&=& \overline{u}_{\infty} , \hspace*{3mm} \partial \Omega_{\bet}
\end{eqnarray*}
In $\Omega_{\bet,\psi}$ we have ;
\begin{equation*}
\left\{
\begin{array}{rcll}
\del u_{\lambda,\beta} &=& \lam \eps ~ \hspace{3mm}&\text{in}~\Omega_{\bet,\psi} \\
u_{\lambda,\beta}&=& \mu~\hspace{3mm} &\text{on } ~ \partial \Omega_{\bet,\psi}
\end{array}
\right.
\end{equation*}
The maximum principle implies
$$\displaystyle \max_{\Omega_{\bet,\psi}}u_{\lambda,\beta}= \displaystyle \max_{\partial \Omega_{\bet,\psi}}u_{\lambda,\beta}= \mu$$
Hence $u_{\lambda,\beta} < \mu$ in $\Omega_{\bet,\psi}$. In $\Omega_{\bet} \setminus \overline{ \Omega_{\bet,\psi}}$ we have ;
\begin{equation*}
\left\{
\begin{array}{rcll}
\del u_{\lambda,\beta} &=& \lam ~ \hspace{3mm}&\text{in}~\Omega_{\bet} \setminus \overline{ \Omega_{\bet,\psi}}\\
u_{\lambda,\beta}&=& \mu~\hspace{3mm} &\text{on } ~ \partial \Omega_{\bet,\psi}\\
u_{\lambda,\beta}&=& \overline{u}_{\infty}\hspace{3mm} &\text{on}~ \partial \Omega_{\bet}
\end{array}
\right.
\end{equation*}
As $\overline{u}_{\infty} > \mu$, then
$$ \displaystyle \min_{\Omega_{\bet} \setminus \overline{ \Omega_{\bet,\psi}}} u_{\lambda,\beta}=\displaystyle \min_{\partial \overline{\Omega_{\bet,\psi}}}u_{\lambda,\beta}=\mu$$
So $u_{\lambda,\beta} > \mu $ in $\Omega_{\bet} \setminus \Omega_{\bet,\psi}$. Therefore, the function $u_{\lambda,\beta}$ satisfies
\begin{equation*}
\left\{
\begin{array}{rcll}
\del u_{\lambda,\beta} &=& \lam ( \eps +(1-\eps ) H(u_{\lambda,\beta}-\mu))~ \hspace{3mm}&\text{in}~\Omega_{\bet} \\
u_{\lambda,\beta}&=& \overline{u}_{\infty}~\hspace{3mm} &\text{on } ~ \partial \Omega_{\bet}
\end{array}
\right.
\end{equation*}
To conclude with the existence of solutions of problem (\ref{sbeta}), we need to prove the existence of function $\psi$ satisfiying the equation
$$u_{\lambda,\beta}(\psi(\theta),\theta)=\mu$$
The variation of the domain $\Omega_{\bet}$ suggests to use an adequate transformation which maps the changing domain into a fixed domain and solves the governing equations in the mapped domain. We consider the transformation $T_{\bet}$ ;
\begin{eqnarray*}
	T_{\bet} : \Omega_{\bet} &\to& \Omega_{0}=B(0,R)\\
	(r,\theta) &\to& (\overline{r},\theta)=(r+r \frac{\bet}{R},\theta)
\end{eqnarray*}
where $(r,\theta)$ is the coordinates in $\Omega_{\bet}$ and $(\overline{r},\theta)$ the coordinates in $\Omega_0$. For a small $\bet$, the transformation $T_{\bet}$ is a diffeomorphism of class $C^2$ of the domain $\Omega_{\bet}$ into $\Omega_0$.
\\The mapping $T_{\bet}$ transform $S_{\bet}$ into $S_0$, hence
$$ T_{\bet}(\psi(\theta),\theta)=(f(\theta),\theta)\hspace*{3mm}\text{where}\hspace{2mm} f \in S_0$$
Using the relation
$$\overline{u}_{\lambda}(T_{\bet}(r,\theta))=u_{\lambda,\beta}(r,\theta)$$
then the problem (\ref{betaxis}) and the equation $u_{\lambda,\beta}(\psi(\theta),\theta)=\mu$ is equivalent to the problem
\begin{equation}\label{0bar}
\left\{
\begin{array}{rcll}
\del \overline{u}_{\lambda}+\delta_{\beta}\overline{u}_{\lambda} &=& \lam(\eps +(1-\eps) \chi_{\Omega_{\bet} \setminus \Omega_{\bet,f}}) ~ \hspace{3mm}&\text{in}~\Omega_{0} \\
\overline{u}_{\lambda}&=& \overline{u}_{\infty}~\hspace{3mm} &\text{on } ~ \partial \Omega_{0}
\end{array}
\right.
\end{equation}
with the equation
\begin{equation}\label{f}
\overline{u}_{\lambda}(f(\theta),\theta)=\mu
\end{equation}
where
\begin{eqnarray*}
\delta_{\beta}&=&\frac{\bet}{R}(2+\frac{\bet}{R})\frac{\partial^2}{\partial \overline{r}^2}+\frac{2 \bet}{\overline{r} R} \frac{\partial}{\partial \overline{r}}\\
&+&\frac{1}{\overline{r}^2} \Big[a_{ij}(\theta)\Big[\frac{r}{R} \frac{\partial \bet}{\partial \theta_j}\Big(\frac{\partial^2}{\partial \overline{r} \partial \theta_i}+\frac{1}{R(1+\frac{\bet}{R})}\frac{\partial \bet}{\partial \theta_i}\frac{\partial}{\partial \overline{r}}+\frac{r}{R}\frac{\partial \bet}{\partial \theta_i}\frac{\partial^2}{\partial \overline{r}^2}\Big)\\
&+&\frac{r}{R} \frac{\partial^2\bet}{\partial \theta_j \partial \theta_i}\frac{\partial}{\partial \overline{r}}+\frac{r}{R}\frac{\partial \bet}{\partial \theta_i}\frac{\partial^2}{\partial \theta_j \partial \overline{r}}\Big]+b_i(\theta)\big[\frac{r}{R}\frac{\partial \bet }{\partial \theta_i}\frac{\partial}{\partial \overline{r}}\big]\Big]
\end{eqnarray*}
with $a_{i,j}$, $b_i$ $\in C^2(S)$, for $i,j=1,2$.\\
Thus, to solve problem (\ref{sbeta}), it is sufficient to prove the existence of function $f$
satisfying (\ref{f}). In fact, using the implicit function theorem, we prove that the equation
$$\overline{u}_{\lambda}(f(\theta),\theta)-\mu=0$$
is uniquely solvable in a given small neighborhood. This is the subject of the following section.
\section{Resolution of integral equation}
For $f \in S_0$, the solution $\overline{u}_{\lambda}$ corresponding to (\ref{0bar}) has an integral representation, given by:
$$\overline{u}_{\lambda}(x)= \overline{u}_{\infty} \int_{\partial \Omega_0} P(x,y) ds(y)+\lam \int_{\Omega_0}\Big(\eps+(1-\eps)\chi_{\Omega_{\bet} \setminus \Omega_{\bet,f}}(y) \Big ) G(x,y) dy-\int_{\Omega_0} \delta_{\beta}\overline{u}_{\lambda}(y) G(x,y)dy$$
where $P$ is the Poisson kernel and $G$ is the Green function. We consider polar coordinates and we define the operator
\begin{eqnarray*}
J :\RR^+ \times S_0 \times \RR^+ \times D &\to&\RR\\
 (\overline{u}_{\infty},f,\mu,\bet) &\to& \overline{u}_{\lambda}(f(\theta),\theta)-\mu
\end{eqnarray*}
where $D$ is the neighborhood of zero in $C(S)$;
\begin{eqnarray*}
J(\overline{u}_{\infty},f,\mu,\bet)(\theta)&=&\overline{u}_{\infty} \int_{\partial \Omega_0} P(f(\theta),\theta,\theta') d \theta'+\lam \int_{\Omega_0}\Big(\eps+(1-\eps)\chi_{\Omega_{\bet} \setminus \Omega_{\bet,f}}(r',\theta') \Big ) G(f(\theta),\theta,r',\theta') dr' d \theta'\\
&-&\int_{\Omega_0} \delta_{\beta}\overline{u}_{\lambda}(r',\theta') G(f(\theta),\theta,r',\theta')dr' d \theta'-\mu\\
&=&\overline{u}_{\infty} \int_{S} P(f(\theta),\theta,\theta') d \theta'+\lam (1-\eps) \int_{S} \int_{f(\theta')}^{R} G(f(\theta),\theta,r',\theta') r'^2 dr' d \theta'\\
&+&\int_{\Omega_0} \big(\lam \eps -\delta_{\beta}\overline{u}_{\lambda}(r',\theta')\big) G(f(\theta),\theta,r',\theta')dr' d \theta'-\mu
\end{eqnarray*}
Then we have the following result.
\begin{lem}\label{lemdiff}
The operator $J$ is continuously differentiable with respect to the second variable.
\end{lem}
\pr
Let $D_j J$ denote the Fréchet derivative of $J$, with respect to the variable of order $j ( j = 1, 2, 3, 4)$. Let $\varphi(\theta)$ be a small perturbation of $f(\theta)$, then We shall prove that $D_2J(\overline{u}_{\infty},f,\mu,\bet)$ is given by
\begin{eqnarray*}
D_2 J(\overline{u}_{\infty},f,\mu,\bet)\varphi(\theta)&=&\overline{u}_{\infty} \int_{S} \frac{\partial P}{\partial r}(f(\theta),\theta,\theta')\varphi(\theta') d \theta'+\lam (1-\eps) \int_{S} \int_{f(\theta')}^{R} \frac{\partial G}{\partial r}(f(\theta),\theta,r',\theta') \varphi(\theta') r'^2 dr' d \theta'\\
&+&\int_{\Omega_0} \big(\lam \eps -\delta_{\beta}\overline{u}_{\lambda}(r',\theta')\big) \frac{\partial G}{\partial r}(f(\theta),\theta,r',\theta') \varphi(\theta') dr' d \theta'\\
&-&\lam (1-\eps) \int_{S} G(f(\theta),\theta,f(\theta'),\theta') \big[f(\theta')\big]^2 \varphi(\theta') d \theta'
\end{eqnarray*}
Let $L = D_2J$. Then $$J(\overline{u}_{\infty},f+\varphi,\mu,\bet)(\theta)-J(\overline{u}_{\infty},f,\mu,\bet)(\theta)-L \varphi(\theta)$$
\begin{eqnarray*}
&=&\overline{u}_{\infty} \int_{S} P(f(\theta)+\varphi(\theta),\theta,\theta') d \theta'-\overline{u}_{\infty} \int_{S} P(f(\theta),\theta,\theta') d \theta'-\overline{u}_{\infty} \int_{S} \frac{\partial P}{\partial r}(f(\theta),\theta,\theta')\varphi(\theta') d \theta'\\
&+&\lam (1-\eps) \int_{S} \int_{f(\theta')+\varphi(\theta')}^{R} G(f(\theta)+\varphi(\theta),\theta,r',\theta') r'^2 dr' d \theta'-\lam (1-\eps) \int_{S} \int_{f(\theta')}^{R} G(f(\theta),\theta,r',\theta') r'^2 dr' d \theta'\\
&-&\lam (1-\eps) \int_{S} \int_{f(\theta')}^{R} \frac{\partial G}{\partial r}(f(\theta),\theta,r',\theta') \varphi(\theta') r'^2 dr' d \theta'\\
&+&\int_{\Omega_0} \big(\lam \eps -\delta_{\beta}\overline{u}_{\lambda}(r',\theta')\big) G(f(\theta)+\varphi(\theta),\theta,r',\theta')dr' d \theta'-\int_{\Omega_0} \big(\lam \eps -\delta_{\beta}\overline{u}_{\lambda}(r',\theta')\big) G(f(\theta),\theta,r',\theta')dr' d \theta'\\
&-&\int_{\Omega_0} \big(\lam \eps -\delta_{\beta}\overline{u}_{\lambda}(r',\theta')\big) \frac{\partial G}{\partial r}(f(\theta),\theta,r',\theta') \varphi(\theta') dr' d \theta'\\
&+&\lam (1-\eps) \int_{S} G(f(\theta),\theta,f(\theta'),\theta') \big[f(\theta')\big]^2 \varphi(\theta') d \theta'
\end{eqnarray*}
We added and subtracted the following term
$$\lam (1-\eps) \int_{S} \int_{f(\theta')}^{R}G(f(\theta)+\varphi(\theta),\theta,r',\theta')r'^2 dr'$$
then we have
$$J(\overline{u}_{\infty},f+\varphi,\mu,\bet)(\theta)-J(\overline{u}_{\infty},f,\mu,\bet)(\theta)-L \varphi(\theta)=I_1+I_2+I_3+I_4$$
where
\begin{eqnarray*}
I_1&=&\overline{u}_{\infty} \int_{S} \Big[ P(f(\theta)+\varphi(\theta),\theta,\theta')- P(f(\theta),\theta,\theta')- \frac{\partial P}{\partial r}(f(\theta),\theta,\theta')\varphi(\theta') \Big]d \theta'\\
I_2&=&\lam (1-\eps) \int_{S} \int_{f(\theta')}^{R} \Big[G(f(\theta)+\varphi(\theta),\theta,r',\theta')-G(f(\theta),\theta,r',\theta')- \frac{\partial G}{\partial r}(f(\theta),\theta,r',\theta') \varphi(\theta')\Big] r'^2 dr' d \theta'\\
I_3&=&\int_{\Omega_0} \Big[G(f(\theta)+\varphi(\theta),\theta,r',\theta')-G(f(\theta),\theta,r',\theta')- \frac{\partial G}{\partial r}(f(\theta),\theta,r',\theta') \varphi(\theta')\Big] \big(\lam \eps -\delta_{\beta}\overline{u}_{\lambda}(r',\theta')\big) dr' d \theta'
\end{eqnarray*}
and
\begin{eqnarray*}
I_4&=&\lam (1-\eps) \int_{S} \Big[-\int_{f(\theta')}^{R}G(f(\theta)+\varphi(\theta),\theta,r',\theta')r'^2 dr'+\int_{(f+\varphi)(\theta')}^{R}G(f(\theta)+\varphi(\theta),\theta,r',\theta')r'^2 dr'\Big] d \theta'\\
&+&\lam (1-\eps) \int_{S} G(f(\theta),\theta,f(\theta'),\theta')\big[ f(\theta')\big]^2 \varphi(\theta') d \theta'\\
&=&-\lam (1-\eps) \int_{S} \Big[ \int_{f(\theta')}^{(f+\varphi)(\theta')}G(f(\theta)+\varphi(\theta),\theta,r',\theta')r'^2 dr'-G(f(\theta),\theta,f(\theta'),\theta')\big[ f(\theta')\big]^2 \varphi(\theta') \Big] d \theta'
\end{eqnarray*}
Using Taylor's theorem, we have $I_1,I_2,I_3,I_4 =o(||\varphi||_{\infty})$ when $||\varphi||_{\infty} \to 0$.
\\Hence,
$$D_2 J(\overline{u}_{\infty},f,\mu,\bet)\varphi(\theta) = \frac{\partial \overline{u}_{\lambda}}{\partial r}(f(\theta),\theta) \varphi(\theta)-\lam (1-\eps) \int_{S} G(f(\theta),\theta,f(\theta'),\theta') \big[f(\theta')\big]^2 \varphi(\theta') d \theta'$$
This proves lemma \ref{lemdiff}\cqd
\begin{lem}\label{leminv}
Suppose that $r_0\neq r_{\lambda_2},$ then the operator $D_2 J$ has a bounded inverse in the neighborhood of $(u_{\infty},r_0,\mu,0)$
\end{lem}
\pr
We know that
$$D_2 J(u_{\infty},r_0,\mu,0)\varphi(\theta) = \frac{\partial u_{\lambda}}{\partial r}(r_0,\theta) \varphi(\theta)-\lam (1-\eps) \int_{S} G(r_0,\theta,r_0,\theta') \big[r_0\big]^2 \varphi(\theta') d \theta'$$
with $u_{\lambda}$ is solution of
\begin{equation*}
\hspace*{0.5cm}
\left\{
\begin{array}{rcll}
\del u &=& \lam ( \eps +(1-\eps ) H(u-\mu))~ \hspace{3mm}&\text{in}~\Omega_0\\
u&=& u_{\infty}~\hspace{3mm} &\text{on } ~ \partial \Omega_0
\end{array}
\right.
\end{equation*}
This implies that
\begin{eqnarray*}
D_2 J(u_{\infty},r_0,\mu,0)\varphi(\theta) &=&  \frac{\lam \eps}{3} r_0  \varphi(\theta)-\lam (1-\eps) \int_{S} G(r_0,\theta,r_0,\theta') \big[r_0\big]^2 \varphi(\theta') d \theta'\\
&=&\lam r_0 \Big(\frac{\eps}{3} I -(1-\eps) r_0 K\Big) \varphi(\theta)
\end{eqnarray*}
where $K$ is the compact operator defined on $C(S)$ by
 $$ K \varphi (\theta)=\int_{S} G(r_0,\theta,r_0,\theta') \varphi(\theta') d \theta'$$
Then we have, $D_2 J(u_{\infty},r_0,\mu,0)$ is invertible provided
\begin{equation}
\frac{\eps}{3} -(1-\eps) r_0 ~\sigma \neq 0 \hspace*{3mm}\text{for any $\sigma$ eigenvalue of $K$}
\end{equation}
We use the expansion of the Green’s function in spherical harmonics (see Appendix. B)
\begin{equation*}
G(r_0,\theta,r_0,\theta')= \frac{1}{r_0} \Big[\frac{r_0-R}{4 \pi R}+ \sum_{l=1}^{\infty} \frac{1}{2l+1} \Big(\Big(\frac{r_0}{R}\Big)^{2l+1}-1\Big)  \sum_{m=1}^{2l+1} Y_{lm}(\theta) Y_{lm}(\theta') \Big]
\end{equation*}
$Y_{lm}$ are the spherical harmonic functions of degree $l$ in dimension $3$. From this expression, we read the eigenvalues $\sigma_l$ of $K$ to be
\begin{equation*}
\sigma_l=\frac{1}{r_0}~ \frac{1}{2l+1} \Big(\Big(\frac{r_0}{R}\Big)^{2l+1}-1\Big) \hspace*{3mm}\text{for $l \in \mathbb{N}$}
\end{equation*}
For $l=0$,
$$ \frac{\eps}{3} -(1-\eps) r_0 ~\sigma_0=\frac{\eps}{3} -(1-\eps) \Big(\frac{r_0}{R}-1\Big)=0$$
This implies $r_0=\frac{2\eps-3}{3(\eps-1)} R=r_{\lam_2}$ which is impossible.
\\For $l \geq 1$, we have
$$\frac{\eps}{3} -(1-\eps) r_0 ~\sigma_l=\frac{\eps}{3}+\frac{\eps-1}{2l+1} \Big(\Big(\frac{r_0}{R}\Big)^{2l+1}-1\Big)$$
we note that $\Big(\frac{r_0}{R}\Big)^{2l+1}-1 <0$ for $r_0 \in (0,R)$, then
\begin{itemize}
\item[*] If $0<\eps \leq 1$, we have $ \Big[(\eps-1) \Big(\Big(\frac{r_0}{R}\Big)^{2l+1}-1\Big)\Big] \geq 0$, so $\Big[\frac{\eps}{3} -(1-\eps) r_0 ~\sigma_l\Big] >0$
\item[*] If $\eps >1$, we have  $ \Big[(\eps-1) \Big(\Big(\frac{r_0}{R}\Big)^{2l+1}-1\Big)\Big] \leq 0$, then
$$\frac{\eps}{3}+\frac{\eps-1}{2l+1} \Big(\Big(\frac{r_0}{R}\Big)^{2l+1}-1\Big) \geq \frac{\eps}{3} + \frac{\eps-1}{3} \Big(\Big(\frac{r_0}{R}\Big)^{2l+1}-1\Big) =\frac{1}{3} +\frac{\eps-1}{3} \Big(\frac{r_0}{R}\Big)^{2l+1} >0$$
\end{itemize}
This proves Lemma \ref{leminv}.\cqd
\begin{rem}
Note that the position of the free boundary $r_0$ is different from $r_{\lambda_2}$ if $\lambda\neq \lambda_2.$ This situation is possible if
\begin{eqnarray*}
(H_2')\hspace*{3mm}\lam &>& \lam_2, \hspace*{3mm}for\hspace*{2mm} \eps \in (3/2,+\infty)\\
(H_3')\hspace*{3mm} \lam &\geq& \lam_1\quad \hbox{and }\quad \lambda\neq \lambda_2   \hspace*{3mm}for\hspace*{2mm} \eps \in (0,3/2].
\end{eqnarray*}
\end{rem}
We have the following result
\begin{theo}
Under assumptions $(H_1)$, $(H_2')$ and $(H_3')$, there exist a neighborhood $V$ of $(u_{\infty},\mu,0)$ in $\RR^+ \times \RR^+ \times C^2(S)$ and a continuous mapping $F : V \to C(S)$ such that
\begin{enumerate}
\item[(i)] $F(u_{\infty},\mu,0)=r_0$
\item [(ii)] $J(\overline{u}_{\infty},F(\overline{u}_{\infty},\mu,\beta),\mu,\beta)=0.$
\end{enumerate}
\end{theo}
\prrr 51
\\Using Lemma \ref{lemdiff} and \ref{leminv} and the implicit function theorem, then we conclude to the
 existence of function $F$ depending on $\overline{u}_{\infty}$, $\mu$ and $\beta$ such that $F(\overline{u}_{\infty},\mu,\beta)$ satisfies $J(\overline{u}_{\infty},F(\overline{u}_{\infty},\mu,\beta),\mu,\beta)=0$
\cqd
When $r_0=r_{\lam_2}=\frac{2\eps-3}{3(\eps-1)} R$ corresponding to $\lambda=\lambda_2,$  the operator $D_2 J(u_{\infty},r_0,\mu,0)$ is not invertible. So, the implicit function theorem fails and a phenomenon of bifurcation appears. We have the following result
\begin{prop}
Assume that $\lambda=\lambda_2$ for $\eps >0$, i.e there exists $\mu^* >0$ such that $$\mu^*=u_{\infty}-\frac{\eps^2 R^2 \big(\frac{4\eps}{3}-\frac{3}{2}\big)}{27(\eps-1)^2} \lambda,\hspace*{3mm}\text{for} ~\lambda >0, \hbox{and}\quad \eps\neq 1.$$
Let $Z=\{\xi \in C(s), \int_{S}\xi(y) dy =0\}$, then there exists
\begin{enumerate}
\item an interval $I=]-\eta,+\eta[$, $\eta >0$.
\item a continuous functions $\varphi : I \to \RR$ and $\psi : I \to \RR$ with $\varphi(0)=\mu^*$ and $\psi(0)=0$.
\item a neighborhood $V$ of $(\mu^*,0)$ in $\RR \times C(S)$ such that for all $s \in I$, the following pair is a solution of $J(u_{\infty},f,\mu,0)=0$ in $V$
$$(\mu,f)=(\varphi(s),r_{\lambda_2}+s \varphi_{00}+s\psi(s))$$
where $\varphi_{00}$ is given constant.
\end{enumerate}
\end{prop}
The proof of proposition 5.1 is based on the following lemmas.
\begin{lem}
	Let $\varphi_{00}=-\frac{1}{4 \pi}$, for $\mu^* >0$, the operator $D_2J(u_{\infty},r_{\lambda_2},\mu^*,0)$ has a one dimensional null space spanned
	by $\varphi_{00}$, while its range has codimension one coinciding with the null space of the continuous linear functional $$\phi(\xi)=\int_S \xi(y) \varphi_{00} dy$$.
\end{lem}
\pr
The operator $D_2J(u_{\infty},r_{\lambda_2},\mu^*,0)$ is not invertible, so
$$D_2J(u_{\infty},r_{\lambda_2},\mu^*,0) \varphi_{00} = \lam r_{\lam_2} \Big(\frac{\eps}{3} I -(1-\eps) r_{\lam_2} K\Big) \varphi_{00} =0$$
This gives that the kernel of $D_2J(u_{\infty},r_{\lambda_2},\mu^*,0)$ is a one dimensional space spanned by $\varphi_{00}$. The function $\varphi_{00}$ is the first eigenfunction corresponding to the eigenvalue $\sigma_0$. Since the operator K is compact, the equation
$$D_2J(u_{\infty},r_{\lambda_2},\mu^*,0) \varphi (\theta)= \xi (\theta),\hspace*{3mm}\text{for}~~ \varphi \in C(S)$$
has a solution if $\xi$ is orthogonal to $\varphi_{00}$. Let
$$\phi(\xi)=\int_S \xi(y) \varphi_{00} dy$$
it becomes apparent that
$$Im ~D_2J(u_{\infty},r_{\lambda_2},\mu^*,0)=Ker~ \phi$$
\cqd
\begin{lem}\cite{Bensid1}
The mixed derivative $D_3 D_2 J(u_{\infty},f,\mu,0)$ exists and is continuous in
a neighborhood of $(r_{\lam_2},\mu^*)$.
\end{lem}
\begin{lem}
$D_3 D_2 J(u_{\infty},r_{\lambda_2},\mu^*,0)\varphi_{00}$ does not belong to the range of $D_2 J(u_{\infty},r_{\lambda_2},\mu^*,0)$
\end{lem}
\pr
When $\beta =0$, we have
$$D_3 D_2 J(u_{\infty},r_{\lambda_2},\mu^*,0)\varphi_{00}=\frac{\partial}{\partial \mu} \Big(\frac{\partial u_{\lam}}{\partial r}(f(\theta),\theta) \varphi_{00}\Big) \mid_{(f=r_{\lam_2}, \mu=\mu^*, \beta=0)} ~~~=\frac{\partial}{\partial \mu} \Big(\frac{\partial u_{\lam}}{\partial r}(r_{\lam_2},\theta) \varphi_{00}\Big)\neq 0$$
this completes the proof.\cqd
To conclude the proof, we take $Z=\{\xi \in C(s), \int_{S}\xi(y) \varphi_{00} dy =0\}$ and we remark that all the hypothesis of bifurcation’s theorem of Crandall-Rabinowitz  are satisfied. See \cite{rabi}.
\section{Appendix}
\subsection*{Appendix A. Cardano's Method}\label{appendixa}
Cardano's method provides a technique for solving the general cubic equation
\begin{equation}\label{car}
a x^3+bx^2+c x +d=0,\hspace*{3mm}a \neq 0
\end{equation}
Cardano's methods involves the following steps:
\begin{enumerate}
\item Eliminate the square term by the substitution $x=z-\frac{b}{3a}$, the method often begins with an equation in the reduced form (the depressed polynomial)
\begin{equation}\label{car2}
z^3+p z+q=0
\end{equation}
where $$p=\frac{-b^2}{3a^2}+\frac{c}{a}, \hspace*{3mm} q=\frac{b}{27 a}\Big( \frac{2b^2}{a^2}-\frac{9c}{a}\Big)+\frac{d}{a}$$
\item Setting $\Delta = q^2+\frac{4}{27}p^3$ the discriminant of the reduced equation, there are three possible solutions to the equation (\ref{car}),
\begin{enumerate}
	\item[(i)] If $\Delta >0$ the equation (\ref{car}) has a real solution, given by
	$$x=\Big(\frac{-q+\sqrt{\Delta}}{2}\Big)^{\frac{1}{3}}+\Big(\frac{-q-\sqrt{\Delta}}{2}\Big)^{\frac{1}{3}}-\frac{b}{3a}$$
	\item [(ii)] If $\Delta = 0$ the equation (\ref{car}) has two real solutions, given by
	$$x_1=\frac{3q}{p}-\frac{b}{3a},\hspace*{3mm}x_2=\frac{-3q}{p}-\frac{b}{3a}$$
	\item[(iii)] If $\Delta < 0$ the equation (\ref{car}) has three real solutions, given by
	$$x_k=2 \sqrt{\frac{-p}{3}} cos \Big(\frac{1}{3} arccos\Big( \frac{-q}{2} \sqrt{\frac{27}{-p^3}}\Big)+\frac{2k\pi}{3}\Big)-\frac{b}{3a}\hspace*{3mm}\text{where $k \in \{0,1,2\}$}$$
\end{enumerate}
\end{enumerate}
Hence, we solve the equation (\ref{12}) using Cardano’s method
\begin{equation*}
\frac{-(\eps-1)}{R} r^3+(\eps-\frac{3}{2}) r^2+ \frac{1}{2} R^2-\frac{3(u_{\infty}-\mu)}{\lam}=0
\end{equation*}
for $r \in (0,R)$, we have
$$a=\frac{-(\eps-1)}{R},~ b=\eps-\frac{3}{2},~c=0,~d=\frac{1}{2} R^2-\frac{3(u_{\infty}-\mu)}{\lam}$$
setting $x=z-\frac{b}{3a}=z+\frac{(\eps-\frac{3}{2})R}{3(\eps-1)}$, then we get the equation (\ref{car2}) with
$$  p=\frac{-(\eps-\frac{3}{2})^2}{3(\eps-1)^2}R^2, \hspace*{3mm}
q=\frac{-2(\eps-\frac{3}{2})^3 R^3}{27(\eps-1)^3}-\frac{ R \big( \frac{1}{2} R^2-\frac{3(u_{\infty}-\mu)}{\lam} \big)}{(\eps-1)}$$
and
\begin{eqnarray*}
\Delta&=&\Big(\frac{-2(\eps-\frac{3}{2})^3 R^3}{27(\eps-1)^3}-\frac{ R \big( \frac{1}{2} R^2-\frac{3(u_{\infty}-\mu)}{\lam} \big)}{(\eps-1)}\Big)^2-\frac{4}{27} \Big(\frac{(\eps-\frac{3}{2})^2}{3(\eps-1)^2}R^2\Big)^3\\
&=&\Big(\frac{R^2(\frac{1}{2}R^2-\frac{3(u_{\infty}-\mu)}{\lam})}{2(\eps-1)^2}\Big)\Big(\frac{1}{2}R^2+\frac{4(\eps-\frac{3}{2})^3}{27(\eps-1)^2}R^2-\frac{3(u_{\infty}-\mu)}{\lam}\Big)
\end{eqnarray*}
Then we study the sign of $\Delta$ according to the values of $\lam_1, \lam_2$ and $\eps$, we get the solutions of the equation (\ref{12}).
\subsection*{Appendix B. Expansion of the Green's function in spherical harmonics}
A fundamental solution for Laplace's equation in spherical coordinates, is
$$G_0(r)=-\frac{1}{4 \pi} r^{-1}$$
Then Green's function for the sphere $B(0,R)$ is
$$G(x,y)=G_0(|x-y|)-G_0\Big(\frac{|x|}{R} ~\Big|\frac{R^2}{|x|^2} x-y\Big|\Big)\hspace*{3mm} \text{see} \cite{dibendetto}$$
On passing to polar coordinates, we find :
\begin{eqnarray*}
G_0(|x-y|)&=&G_0(\big|r_0 \theta-r_0\theta'\big|)= G_0\Big(\sqrt{2r_0^2-2r_0^2 \theta.\theta'}\Big)
	=G_0\Big(r_0\sqrt{2(1-cos \gamma)}\Big)\\
G_0\Big(\frac{|x|}{R} ~\Big|\frac{R^2}{|x|^2} x-y\Big|\Big)&=&G_0\Big(\sqrt{\frac{r_0^2}{R^2}\Big(\frac{R^4}{r_0^4} r_0^2+r_0^2-2\frac{R^2}{r_0^2}r_0^2 \theta.\theta'\Big)}\Big)
= G_0 \Big(R \sqrt{1+\Big(\frac{r_0}{R}\Big)^4-2 \Big(\frac{r_0}{R}\Big)^2 cos \gamma }\Big)	
\end{eqnarray*}
where $cos \gamma=\theta.\theta'$, then
\begin{equation*}
(B_0)\hspace*{3mm} G(r_0,\theta,r_0,\theta')=G_0\Big(r_0\sqrt{2(1-cos \gamma)}\Big)-G_0 \Big(R \sqrt{1+\Big(\frac{r_0}{R}\Big)^4-2 \Big(\frac{r_0}{R}\Big)^2 cos \gamma }\Big)
\end{equation*}
Now expand the second term on the right in Legendre polynomials for $n=3$ (\cite{harmonics} Lemma 18)
\begin{equation*}
(B_1)\hspace*{3mm}	G_0 \Big(R \sqrt{1+\Big(\frac{r_0}{R}\Big)^4-2 \Big(\frac{r_0}{R}\Big)^2 cos \gamma }\Big)=\frac{-1}{4 \pi R} \sum_{l=0}^{\infty} \Big(\frac{r_0}{R}\Big)^{2l} P_l(3,cos\gamma)
\end{equation*}
and apply the Addition Theorem for spherical harmonics, we find
$$ (B_2)\hspace*{3mm} P_l(3,cos\gamma)=\frac{4 \pi}{2l+1} \sum_{m=1}^{2l+1} Y_{lm}(\theta) Y_{lm}(\theta')$$
with the $Y_{lm}$ an orthonormal basis for spherical harmonics of degree $l$, in dimension $3$.
\\Using $(B_1)$ and $(B_2)$, we have
\begin{eqnarray*}
(B_3)\hspace*{3mm}	G_0 \Big(R \sqrt{1+\Big(\frac{r_0}{R}\Big)^4-2 \Big(\frac{r_0}{R}\Big)^2 cos \gamma }\Big) &=& -\frac{1}{R}~ \sum_{l=0}^{\infty}\frac{1}{2l+1} \Big(\frac{r_0}{R}\Big)^{2l} \sum_{m=1}^{2l+1} Y_{lm}(\theta) Y_{lm}(\theta')\\
&=&-\frac{1}{4 \pi R}-\frac{1}{R}~\sum_{l=1}^{\infty}\frac{1}{2l+1} \Big(\frac{r_0}{R}\Big)^{2l}  \sum_{m=1}^{2l+1} Y_{lm}(\theta) Y_{lm}(\theta')
\end{eqnarray*}
For the first term on the right in $(B_0)$ we use Lemma 18 and the Addition Theorem in \cite{harmonics} again, we obtaining
\begin{equation*}
	(B_4)\hspace*{3mm} G_0\Big(r_0\sqrt{2(1-cos \gamma)}\Big)=-\frac{1}{4 \pi r_0}-\frac{1}{r_0}~\sum_{l=1}^{\infty}\frac{1}{2l+1} \sum_{m=1}^{2l+1} Y_{lm}(\theta) Y_{lm}(\theta')
\end{equation*}
The series converge uniformly in $|\theta-\theta'|\geq \eta >0$. Combining $(B_3)$ and $(B_4)$, gives finally expansion of our kernel in spherical harmonics in dimension $3$.
\begin{equation*}
G(r_0,\theta,r_0,\theta')= \frac{1}{r_0} \Big[\frac{r_0-R}{4 \pi R}+\sum_{l=1}^{\infty}\frac{1}{2l+1} \Big(\Big(\frac{r_0}{R}\Big)^{2l+1}-1\Big)  \sum_{m=1}^{2l+1} Y_{lm}(\theta) Y_{lm}(\theta') \Big]
\end{equation*}

\end{document}